\documentclass[reqno]{amsart}
\usepackage{amssymb}
\usepackage{graphicx}
\usepackage{url}
\def\qmod#1#2{{\hbox{}^{\displaystyle{#1}}}\!\big/\!\hbox{}_{
\displaystyle{#2}}}

\def\resto#1#2{{
#1\hskip 0.4ex\vline_{\hskip 0.4ex\raisebox{-0,5ex}
{{${\scriptstyle #2}$}}}}}
\def\at#1#2{
#1\hskip 0.25ex\vline_{\hskip 0.25ex\raisebox{-1.5ex}
{{$\scriptstyle#2$}}}}

\newfam\msbfam

\font\tenmsb=msbm10

\font\sevenmsb=msbm10 at 7pt
\font\fivemsb=msbm10 at 5pt

\textfont\msbfam=\tenmsb
\scriptfont\msbfam=\sevenmsb
\scriptscriptfont\msbfam=\fivemsb

\def\C{{\mathbb C}}

\def\H{{\mathbb H}}

\def\P{{\mathbb P}}

\def\R{{\mathbb R}}

\def\Z{{\mathbb Z}}

\def\union{\mathop{\bigcup}}
\def\qed {\hfill\vrule height6pt width6pt depth0pt \bigskip}

\def\map{\longrightarrow}
\def\textmap#1{\mathop{\vbox{\ialign{
                                 ##\crcr
     ${\scriptstyle\hfil\;\;#1\;\;\hfil}$\crcr
     \noalign{\kern 1pt\nointerlineskip}
     \rightarrowfill\crcr}}\;}}

\newcommand{\cal}{\mathcal}
\def\textlmap#1{\mathop{\vbox{\ialign{
                                 ##\crcr
     ${\scriptstyle\hfil\;\;#1\;\;\hfil}$\crcr
     \noalign{\kern-1pt\nointerlineskip}
     \leftarrowfill\crcr}}\;}}

\newfam\meuffam

\font\tenmeuf=eufm10
\font\sevenmeuf=eufm7
\font\fivemeuf=eufm5

\textfont\meuffam=\tenmeuf
\scriptfont\meuffam=\sevenmeuf
\scriptscriptfont\meuffam=\fivemeuf

\def\germ{\fam\meuffam\tenmeuf}

\def\g{{\germ g}}

\newtheorem{sz}{Satz}[section]
\newtheorem{thry}[sz]{Theorem}
\newtheorem{pr}[sz]{Proposition}
\newtheorem{re}[sz]{Remark}
\newtheorem{co}[sz]{Corollary}
\newtheorem{dt}[sz]{Definition}
\newtheorem{lm}[sz]{Lemma}

\begin{document}
\def\Pr{{\rm Pr}}
\def\tr{{\rm Tr}}
\def\End{{\rm End}}
\def\Aut{{\rm Aut}}
\def\Spin{{\rm Spin}}
\def\U{{\rm U}}
\def\SU{{\rm SU}}
\def\SO{{\rm SO}}
\def\PU{{\rm PU}}
\def\GL{{\rm GL}}
\def\spin{{\rm spin}}
\def\su{{\rm su}}
\def\so{{\rm so}}
\def\ub{\underbar}
\def\pu{{\rm pu}}
\def\Pic{{\rm Pic}}
\def\Iso{{\rm Iso}}
\def\NS{{\rm NS}}
\def\deg{{\rm deg}}
\def\Hom{{\rm Hom}}
\def\Aut{{\rm Aut}}
\def\h{{\germ h}}
\def\Herm{{\rm Herm}}
\def\Vol{{\rm Vol}}
\def\pf{{\bf Proof: }}
\def\id{{\rm id}}
\def\i{{\germ i}}
\def\im{{\rm im}}
\def\rk{{\rm rk}}
\def\ad{{\rm ad}}
\def\h{{\bf H}}
\def\coker{{\rm coker}}
\def\dbar{\bar{\partial}}
\def\Lo{{\Lambda_g}}
\def\niq{=\kern-.18cm /\kern.08cm}
\def\Ad{{\rm Ad}}
\def\RSU{\R SU}
\def\ad{{\rm ad}}
\def\dva{\bar\partial_A}
\def\da{\partial_A}
\def\p{\partial\bar\partial}
\def\sp{\Sigma^{+}}
\def\sm{\Sigma^{-}}
\def\spm{\Sigma^{\pm}}
\def\smp{\Sigma^{\mp}}
\def\Tors{{\rm Tors}}
\def\st{{\rm st}}
\def\s{{\rm s}}
\def\oo{{\scriptstyle{\cal O}}}
\def\ooo{{\scriptscriptstyle{\cal O}}}
\def\sw{Seiberg-Witten }
\def\pa{\partial_A\bar\partial_A}
\def\Dr{{\raisebox{0.15ex}{$\not$}}{\hskip -1pt {D}}}
\def\gr{{\scriptscriptstyle|}\hskip -4pt{\g}}
\def\subsetint{{\  {\subset}\hskip -2.45mm{\raisebox{.28ex}
{$\scriptscriptstyle\subset$}}\ }}
\def\ra{\rightarrow}
\def\pst{{\rm pst}}
\def\sst{{\rm sst}}

\def\kod{{\rm kod}}
\def\degmax{{\rm degmax}}
\def\red{{\rm red}}
\def\ASD{{\rm ASD}}

\title{Families of holomorphic bundles }

\begin{abstract}
The first goal of the article is to  solve several fundamental problems in the theory of holomorphic
bundles over non-algebraic manifolds: For instance we prove that stability and semi-stability are Zariski open properties in families when the Gauduchon degree map is a topological invariant, or when the parameter manifold is compact. Second we show that, for a generically stable family of bundles over a K\"ahler manifold,  the Petersson-Weil form extends as a closed positive current on the whole parameter space of the family. This extension theorem uses classical tools from Yang-Mills theory (e.g. the Donaldson functional on the space Hermitian metrics and its properties). We apply these results to  study  families of bundles over a K\"ahlerian manifold $Y$ parameterized by a non-K\"ahlerian surface $X$, proving that  such families must satisfy very restrictive conditions.  These results play an important role in our program to prove existence of curves on class VII surfaces \cite{Te1}, \cite{Te2}, \cite{Te3}.
\end{abstract}
\author{Andrei Teleman}
\date{\today}
\maketitle

\tableofcontents

\section{Summary of results}

 There are many important results in the theory of holomorphic bundles which suggest that, in  general 
 \\
 \\
{\bf P:}  {\it A moduli space of bundles on a compact complex manifold  inherits fundamental geometric properties from the base manifold.}
 \\  
 
We mention several well-known celebrated results which illustrate this principle. 
\begin{enumerate}
\item  Any moduli space of stable bundles on a projective  manifold has a natural projective  compactification, so it is quasi-projective.\\
This follows from the Gieseker-Maruyama stability theory.
\item Any moduli space of stable bundles on a  K\"ahlerian manifold is K\"ahlerian.\\
This statement is an easy consequence of the Kobayashi-Hitchin correspondence  on K\"ahler manifolds obtained by Donaldson and Uhlenbeck-Yau. 

\item Any moduli space of stable sheaves on a holomorphically symplectic surface is smooth and holomorphically symplectic, by a result of Mukai.
\end{enumerate}

These statements can be reformulated in a different way, as {\it incompatibility theorems}. For instance by the first result one cannot expect to have a non-algebraic compact subvariety in a moduli space of stable bundles over a projective algebraic   manifold. Therefore, there do not exist non-degenerate families of holomorphic bundles over an algebraic base parameterized by a compact non-algebraic parameter manifold.  But one can also ask the other way round: 
\\ 
\\
{\bf Q:}  {\it Let ${\cal M}^\st$ be a moduli space of stable bundles over  
\begin{enumerate}
\item a compact non-algebraic complex  manifold $X$, or
\item  a Gauduchon surface $(X,g)$   with $b_1(X)$ odd. 
\end{enumerate}
Does there exist a compact  positive dimensional  algebraic  variety  $Y$ (respectively a compact positive dimensional  smooth  K\"ahler manifold $(Y,g)$) with a holomorphic embedding $Y\hookrightarrow {\cal M}^\st$? If yes, classify the corresponding bundle families.}
\\
\\
The aim of this section is to explain why answering such  questions is important,   which methods will be used to deal with these problems, and what are the difficulties. Then we will list our results.\\

Our proof of the ``global spherical conjecture" in the case $b_2=1$ \cite{Te1} uses essentially a result of this type, namely: 
\begin{thry}\label{curves} Let $X$ be a compact complex manifold of algebraic dimension 0, $Y$ a Riemann surface,   and let  ${\cal E}$ be a holomorphic rank 2 bundle on  $Y\times X$. Then there exists a reflexive coherent sheaf ${\cal T}_0$ on $X$ of rank 1 or 2, a non-empty Zariski open set $U\subset X$, and for every $y\in Y$ a morphism $e_y:{\cal T}_0\to {\cal E}_y$ which is a bundle embedding (i.e. fiberwise injective) on $U$. 
\end{thry}

Therefore, there are very few families which can occur. For instance, one can have a family of extensions 
$$0\to {\cal L}_0\map {\cal E}\map {\cal M}_0\map 0
$$
(with fixed ${\cal L}_0,\ {\cal M}_0\in\Pic(X)$) parameterized by a curve in the projective space $\P({\rm Ext}^1({\cal M}_0,{\cal L}_0))$. Another typical class of examples (which correspond to the case $\rk({\cal T}_0)=2$) is a one-parameter family of ``elementary transformations", i.e. of kernels of the form
$$0\map {\cal E}_y\map {\cal F}_0\textmap{u_y} (i_{D_0})_*{\cal F}_y\map 0\ , 
$$
where ${\cal F}_0$ is fixed rank 2 bundle on $X$, and $u_y:{\cal F}_0\to(i_{D_0})_*{\cal F}_y$ are epimorphisms onto the torsion sheaves associated with a family of  bundles ${\cal F}_y$ (of rank 1 or 2) defined on a fixed effective divisor $D_0\subset X$. Both methods yield families of bundles with fixed determinant line bundle. \\

Our method to prove existence of curves on class VII surfaces can be extended to higher $b_2$  (see \cite{Te3}) but even for $b_2=2$  important new difficulties arise. For instance, Theorem \ref{curves} above will no longer be sufficient. Instead one has to answer the second question {\bf Q}(2) concerning the existence of compact K\"ahlerian submanifolds embedded in moduli spaces  of stable bundles over non-K\"ahlerian surfaces.
 
Our strategy to answer such questions starts with a very natural idea (see \cite{Te1}):  Suppose that $Y$ has a good complex geometric property ${\cal P}$ (e.g. algebraic or K\"ahlerian\footnote{We agree to call {\it K\"ahlerian} a complex manifold which admits K\"ahler metrics, whereas a {\it K\"ahler} manifold is a complex manifold {\it endowed} with a K\"ahler metric.}) which  $X$ does not have. Given a family ${\cal E}\to Y\times X$ of bundles on $X$ parameterized by  $Y$,  {\it switch the roles of the two terms}, i.e. regard ${\cal E}$  as a family of bundles on $Y$ parameterized by $X$. If the obtained bundles ${\cal E}^x$, $x\in X$ were all stable, then  we would  get a regular holomorphic map $f:X\to {\cal M}^\st$ in a moduli space of bundles on $Y$; the general principle {\bf P}  applies, showing that $f$ must  factorize through a variety  $X_0$ with the property ${\cal P}$. This reduces the problem  to the classification of families on $X_0\times Y$, where $X_0$ is a ``quotient of $X$" having the property ${\cal P}$. Unfortunately, there is no reason to expect the bundles ${\cal E}^x$ to be all stable. 

The interesting  and difficult case   is the {\it generically stable case}, i.e. the case when these bundles are stable {\it for generic} $x\in X$.  Suppose that we are in this case, and put
$$X^\st:=\{x\in X|\ {\cal E}^x\hbox{ is stable }\}\ .
$$
 When ${\cal P}$={\it projective algebraic} and $X$ is a smooth non-algebraic surface, there is an natural way to make our method work: a generically stable  family defines a meromorphic map $f:X\dasharrow \bar {\cal M}^\st$ in the Gieseker compactification of the moduli space, so a regular morphism $\tilde f:\tilde X\to \bar {\cal M}^\st$ defined on a modification $\tilde X$ of $X$, which will also be a non-algebraic surface. Therefore $\tilde f$ must be degenerate. For instance, if $X$ was supposed to have vanishing algebraic dimension, it will follow that $f$ must be constant on $X^\st$.  Therefore the initial family ${\cal E}$ must be generically constant with respect to $x\in X$; it is easy to see that this implies the conclusion of Theorem \ref{curves} with $ {\cal T}_0={\cal O}_X^{\oplus 2}\otimes {\cal U}$ for a line bundle ${\cal U}$. This simple argument  gives a new proof and a generalization of this theorem in the generically stable case. The case of a family which is not generically stable can be handled as in \cite{Te1} using  Lemma \ref{nonstable} (compare to the proof of Theorem \ref{incompatibility} in section \ref{incompatsection}). Therefore Theorem \ref{curves} holds for an arbitrary projective algebraic manifold $Y$.\\\
 
 When ${\cal P}$={\it K\"ahlerian}  and $X$ is a non-K\"ahlerian surface, one could try the same strategy. For a K\"ahler manifold $(Y,g)$ one has an induced K\"ahler structure on any moduli space ${\cal M}^\st$ of stable bundles on $Y$.   In a presence of a family ${\cal E}\to Y\times X$ one can apply the same method as above and  hope to either obtain  an induced singular  K\"ahler structure on $X$ (induced from   ${\cal M}^\st$) or to prove that the family is  degenerate at every point. This method  encounters serious difficulties (see section \ref{remarks} ``remarks and open problems" for details):
 \begin{enumerate}
 \item In the non-algebraic framework it is well-known  that $X^\st$ is open \cite{LT1}  but, to our knowledge,  not that it is {\it Zariski open}.
 \item In the K\"ahlerian non-algebraic framework we do not have a {\it complex geometric} compactification of the moduli space of stable bundles ${\cal M}^\st$. Even for $\dim(X)=2$ it is not known whether the Uhlenbeck compactification   (which is the compactification suggested by the Kobayashi-Hitchin correspondence) is a complex space. Moreover, if one tries to compactify using moduli spaces of semistable sheaves, it is very difficult to prove compactness, because the gauge theoretical  topology on the space of virtual instantons does not compare easily to the topology defined by holomorphic deformation theory for sheaves.  
 \item Even if a complex geometric compactification of ${\cal M}^\st$ was constructed (e.g. by endowing the Uhlenbeck compactification with a complex space structure, or proving compactness of a  certain moduli space of semistable sheaves), it is very hard to see whether the $L^2$-K\"ahler metric on ${\cal M}^\st$ (also called the Petersson-Weil metric) extends to the compactification.
\end{enumerate}

Note that in the non-K\"ahlerian framework   these problems have all  {\it negative} answers:  
\begin{re} For a non-K\"ahler Gauduchon manifold $(Y,g)$ \begin{enumerate}
\item   Stability  with respect to $g$ is open but in general not Zariski open in families. Semistability is not open in general, even with respect to the classical topology.
\item  A moduli space ${\cal M}^\st$ of $PSL(2,\C)$ bundles  cannot be compactified as a compact complex space in general. There exists compact moduli spaces of polystable bundles over  class VII surfaces which identify to compact disks \cite{Sch} or the sphere $S^4$ \cite{Te3}.
\end{enumerate}
\end{re}
This shows that, in bundle theory, many standard algebraic geometric results cannot be automatically extended to the non-algebraic framework; generalizing such standard results requires new proofs adapted to the non-algebraic framework.\\

We will begin the article  by solving the first problem (see section \ref{opensection}).  We will need a much weaker hypothesis than K\"ahlerianity:
\begin{thry}\label{open}
Let $(Y,g)$ be a compact connected Gauduchon manifold, $X$ an arbitrary complex manifold $X$, and ${\cal E}$ a holomorphic rank $r$ bundle on $X\times Y$. Set
$$X^\st:=\{x\in X|\ {\cal E}_x\hbox{ is stable }\}\ ,\ X^\sst:=\{x\in X|\ {\cal E}_x\hbox{ is semistable }\}\ .
$$
 Then
 \begin{enumerate}
 \item $X^\st$ is   open with respect to the classical topology.
 \item   Suppose that   one of the following conditions is satisfied:
 \begin{enumerate}
 \item \label{ti} $\deg_g:\Pic(Y)\to\R$ is a topological invariant (i.e. it vanishes on $\Pic^0(Y)$).
 \item \label{cb} the parameter manifold $X$ is compact.
 \end{enumerate}
Then $X^\st$ and $X^\sst$ are both Zariski open.  
\end{enumerate}
 \end{thry}

Second, in section \ref{PW} we prove a general extensibility theorem, which will allow us to avoid the construction of  a complex geometric K\"ahler compactification of a moduli space of stable bundles over a K\"ahler manifold (hence to avoid the second and third difficulties in the above list). Our result is:

  \begin{thry}\label{extension} Let ${\cal E}$ be a family over $X\times Y$ with $(Y,g)$ K\"ahler. Suppose that the Zariski open set $X^\st$ is non-empty, i.e. ${\cal E}$ is generically stable. The Petersson-Weil form   $pw(\resto{{\cal E}}{X^\st\times Y})$ extends in a  natural way  as a closed positive current on $X$. 
  \end{thry}

The proof will make use of  fundamental  tools  in the Yang-Mills theory which have been developed  by Donaldson  \cite{Do1}.    The idea is to introduce a parameter dependent version $m(K,H)$ of Donaldson's functional on the space of Hermitian metrics on the family ${\cal E}$, and to prove that the function $m(K,H_{he})\in{\cal C}^\infty(X^\st,\R)$ associated with the metric  $H_{he}$ obtained by solving fiberwise the Hermite-Einstein equation on $\resto{{\cal E}}{X^\st\times Y}$  extends as an almost plurisubharmonic function on $X$.

As we explain in section \ref{remarks}, our extensibility result raises many interesting questions. For instance, one can ask: what are the Lelong numbers of this current, or  describe the complex subspace defined by the Nadel's ideal sheaf of $pw({\cal E})$ in terms of the complex geometry of the family ${\cal E}$.  \\
 
Using these preparations, we will be able to prove the following ``incompatibility theorem" for generically stable families of bundles on a K\"ahler manifold parameterized by a a non-K\"ahlerian surface.
\begin{thry}\label{incompatibility} Let $(Y,g)$ be a compact K\"ahler manifold, $X$ a surface with $b_1(X)$ odd and ${\cal E}\to X\times Y$ a generically stable family of bundles on $Y$  parameterized by $X$. Then
\begin{enumerate}
\item The family is degenerate at every point, i.e. the induced map $X^\st\to {\cal M}^{st}$ is either constant or has generically  rank 1. 
\item In the latter case $X$ contains infinitely many compact curves so $a(X)\geq 1$.
\end{enumerate}
\end{thry}
Using this result we will prove the following K\"ahlerian version  of Theorem \ref{curves}.
\begin{co}\label{goal} Let $(Y,g)$ be compact K\"ahler manifold, $X$ a surface with $b_1(X)$ odd and $a(X)=0$. Let ${\cal E}\to X\times Y$ be an arbitrary family of rank 2 bundles on $Y$  parameterized by $X$. Then there exist a locally free  sheaf ${\cal T}_0$  on $X$ of rank 1 or 2, a non-empty open Zariski open set $U\subset X$ and, for every $y\in Y$, a morphism $e_y:{\cal T}_0\to {\cal E}^y$ which is a bundle embedding on $U$.
\end{co}
This result  plays an important role in our recent work about existence of curves on class VII surfaces with $b_2=2$ \cite{Te3}. 
In the proof of Theorem \ref{incompatibility} and Corollary \ref{goal} we will need a Demailly type self-intersection inequality for closed positive currents on complex surfaces. It is just a version of Demailly's self-intersection inequality (see Corollary 7.6 in \cite{De1}), which holds in arbitrary dimension but requires K\"ahlerianity. Our non-K\"ahlerian version -- valid only on surfaces -- cannot be found in the literature, so we wrote an appendix dedicated to this result. \\

This article has been conceived as a part of our program to prove existence of curves on class VII surfaces. In order to proceed with this program and to extend it to surfaces with $b_2=2$ we had to overcome new difficulties,  and to prove theorems  concerning the general theory of holomorphic bundles.  Therefore we believe that many of these results are of independent interest.   
\\ \\
{\bf Acknowledgements:} The author is indebted to Dan Popovici, Tien-Cuong Dinh and S\'ebastien Boucksom  for their valuable suggestions and comments on  the complex-analytic part of the article. Dan Popovici  suggested a short proof -- based on recent results of S\'ebastien Boucksom -- of  the self-intersection inequality given in Corollary  \ref{dan}.  From Tien-Cuong Dinh I learnt several extension theorems for plurisubharmonic functions; one of these plays a crucial role in   the proof of Theorem \ref{extension}. I also benefited from very useful discussions about moduli spaces of bundles  with Misha Verbitsky, Matei Toma and Nicholas Buchdahl, and from an interesting exchange of e-mails with Jun Li and Gang Tian concerning the compactification problem  and the extensibility of the $L^2$-metric in the non-algebraic framework.

\section{Openness and Zariski openness properties of stability} 
\label{opensection}

\subsection{Properness theorems}\label{max}

Let $Y$ be compact complex $n$-dimensional manifold and ${\cal E}$ a holomorphic rank $r$ vector bundle over $Y$.  The connected components $\Pic^c(Y)$, $c\in NS(Y)$ of   the Picard group $\Pic(Y)$ might be non-compact in general, so an analytic subset of such a component is not necessary compact.
We denote by ${\cal D}ou(Y)$ the Douady space of all effective divisors of $Y$ \cite{Dou}. The  natural map  
$$n_Y:{\cal D}ou(Y)\to \Pic(Y)\ ,\ D\mapsto [{\cal O}(D)]$$
 is a projective fibration over its image. It is well known, that, 
 \begin{pr}\label{compactdouady} For any Hermitian metric  $g$ on $Y$, the topological subspaces 
$${\cal D}ou(Y)_{\leq s}:=\{D\in {\cal D}ou(Y)|\ vol_g(D)\leq s\}
$$
 of ${\cal D}ou(Y)$ are compact.   %
 \end{pr}
  For non-K\"ahlerian metrics, these sets are not complex   analytic in general, because the real function $vol_G$ might   be non-constant on the connected components of ${\cal D}ou(Y)$.   Proposition \ref{compact} implies that $n_Y({\cal D}ou(Y)_{\leq s})$ intersects only finitely many  components of $\Pic(Y)$ and each intersection is compact (but in general not analytic!).
\\

\begin{dt} Let ${\cal E}$ be a holomorphic bundle on $Y$. The {\it Brill-Noether locus associated with} ${\cal E}$ is the analytic subset of $\Pic(Y)$ defined by
$$BN({\cal E}):=\{[{\cal L}]\in\Pic(Y)|\  H^0({\cal L}^\vee\otimes{\cal E})\ne 0\}\ .
$$
\end{dt}
In general,  for a fixed $c\in NS(Y)$, the intersection $BN^c({\cal E}):=BN({\cal E})\cap \Pic^c(Y)$ might be non-compact.  \\
 
 Denote by ${\cal Q}uot^1_{\rm lf}({\cal E})$ the open subspace of the Douady Quot space ${\cal Q}uot({\cal E})$ consisting of quotients with rank 1 locally free kernel.  ${\cal Q}uot^1_{\rm lf}({\cal E})$ comes with a natural holomorphic map 
$$\kappa:{\cal Q}uot^1_{\rm lf}({\cal E})\to \Pic(Y)$$
 which assigns to a quotient its kernel.  The image of this map is just   the Brill-Noether locus $BN({\cal E})$, and the fiber over $[{\cal L}]\in BN({\cal E})$ (endowed with the reduced structure) is just $\P(H^0({\cal L}^\vee\otimes{\cal E}))$. In other words, ignoring possible non-reduced structures, ${\cal Q}uot^1_{\rm lf}({\cal E})$ can be regarded as a projective fibration over $BN({\cal E})$.\\

Put $Z:=\P({\cal E}^\vee)$ and let $q:Z\to Y$ be the natural projection. Let ${\cal L}$ be a line bundle and $\varphi\in H^0({\cal L}^\vee\otimes{\cal E})\setminus\{0\}$.   Using the  natural isomorphisms 
\begin{equation}
H^0({\cal E}\otimes{\cal L}^\vee)=H^0(q_*({\cal O}_Z(1)\otimes q^*({\cal L}^\vee)))=H^0({\cal O}_Z(1)\otimes q^*({\cal L}^\vee))\ ,\end{equation}
we can identify $\varphi$ with a non-trivial section $\tilde\varphi$ in the line bundle ${\cal O}_Z(1)\otimes q^*({\cal L}^\vee)$ over $Z$. Denote by $D(\varphi)\in  {\cal D}ou(Z)$   the corresponding divisor. 
The embedding  $a:\Pic(Y)\to\Pic(Z)$ given by  
$$[{\cal L}]\mapsto[{\cal O}_Z(1)\otimes q^*({\cal L}^\vee)]\ . $$
is closed and open (it identifies $\Pic(Y)$ with a union of connected components of $\Pic(Z)$). 
\begin{pr} \label{iso} The map $\varphi\mapsto D(\varphi)$ defines an isomorphism
$$F:{\cal Q}uot^1_{\rm lf}({\cal E})\textmap{\simeq} n_Z^{-1}(a(\Pic(Y))\subset{\cal D}ou(Z)\ .
$$
 \end{pr} 
\pf  Put ${\cal D}:=n_Z^{-1}(a(\Pic(Y))\subset{\cal D}ou(Z)$ to save on notations\\ \\
Step 1. First of all note that $F$ is obviously a bijection.\\ \\
Step 2. Both spaces ${\cal D}$ and ${\cal Q}uot^1_{\rm lf}({\cal E})$ represent  functors defined on the category of complex spaces. The first space represents the functor $A$ which assigns to a complex space $S$ the set of divisors $D\subset S\times Z$ such that ${\cal O}_D$ is flat over $S$ and such that ${\cal O}(D_s)\in a(\Pic(Y))$ for all $s\in S$. The second space represents the functor $B$ which assigns to $S$ the set of quotients of ${\cal E}_S:=p_S^*({\cal E})$ which are flat over $S$ and have invertible kernel  over $S\times Y$.\\ \\
Step 3. The two functors $A$ and $B$ are  isomorphic: Indeed,  we put $q_S:=\id_S\times q$ and we apply Step 1 to the pair $(S\times Z\textmap{ q_S}S\times Y,{\cal E}_S)$ instead of $(Y,{\cal E})$. For a line bundle ${\cal L}$ over $S\times Y$ we get again an isomorphism
$$\Hom({\cal L},{\cal E}_S)=H^0({\cal L}^\vee\otimes {\cal E}_S)=H^0( q^*_S({\cal L}^\vee)(1))\simeq \Hom({\cal O}_{S\times Z},   q^*({\cal L}^\vee)(1))\ .
$$
 
 On non-reduced or reducible complex spaces non-trivial morphisms ${\cal L}\to {\cal E}$ are not necessary  sheaf monomorphisms. But it is easy to see that, via the above isomorphism, sheaf monomorphisms ${\cal L}\to{\cal E}$ correspond to sheaf monomorphisms ${\cal O}_{S\times Z}\to  q^*({\cal L}^\vee)(1)$.
 
 It remains to check that the two flatness conditions are equivalent. Let ${\cal Q}$ be the quotient of a monomorphism ${\cal L}\to {\cal E}_S$, and ${\cal Q}'$ the quotient of the corresponding monomorphism ${\cal O}_{S\times Z}\to q^*({\cal L}^\vee)(1)$. By the local flatness criterion (\cite{Fi} Theorem 3.14) the $S$-flatness of ${\cal Q}$ at $(s,y)$ is equivalent to the condition ${\rm Tor}_1^{{\cal O}_s}(\C_s,{\cal Q}_{s,y})=0$. But, since the projection $S\times Y\to S$ is flat, we have 
 $${\rm Tor}_1^{{\cal O}_s}(\C_s,{\cal Q}_{(s,y)})={\rm Tor}_1^{{\cal O}_{(s,y)}}(\C_s\otimes_{{\cal O}_s}{\cal O}_{(s,y)},{\cal Q}_{(s,y)})\ ,$$
so ${\cal Q}$ is flat over $S$ if and only if  the sheaf ${\cal T}or_1({\cal O}_{\{s\}\times Y},{\cal Q})$ vanishes for every $s\in S$. This is equivalent to the injectivity of  the induced morphism 
 $${\cal L}_s:={\cal L}_{\{s\}\times Y}\to [{\cal E}_S]_{\{s\}\times Y}\simeq {\cal E}$$
over $Y$,  for every $s\in S$. Similarly, the flatness of ${\cal Q}'$ is equivalent to the injectivity of the morphism
 $${\cal O}_Z\simeq {\cal O}_{\{s\}\times Z}\to q^*({\cal L}_s^\vee)(1)
 $$
 over $Z$, for every $s\in S$. It suffices to notice that the two injectivity conditions are equivalent.
\qed

 Choose a Gauduchon metric $g$ on $Y$. We endow $Z=\P({\cal E}^\vee)$ with a {\it Hermitian} metric in the following way. We choose a Hermitian metric $h$ on ${\cal E}$ and we denote by $A_h$ the corresponding Chern connection. The connection $A_h$ defines a connection $B_h$ in the locally trivial bundle $Z=\P({\cal E}^\vee)\to Y$. The $B_h$-horizontal spaces are invariant under the almost complex structure of $Z$. The vertical tangent bundle $T^V_Z$ comes with a natural Hermitian structure induced by the Fubini-Study metric on the projective fibers; we normalize this metric such that the volume of these fibers is 1, and we denote by  $\omega_{FS}$ the corresponding $(1,1)$-form on this bundle. Put
 $$\Omega_{g,h}:=q^*(\omega_g)+[q^V_h]^*(\omega_{FS})\ ,
 $$
  where $q^V_h$ stands for the projection on the vertical bundle associated with the connection $B_h$. $\Omega_{g,h}$ is the form of a Hermitian metric  $H_{g,h}$ on  $Z$ with respect to which the vertical and the $B_h$-horizontal tangent bundles become orthogonal. 
  We do not claim that $H_{g,h}$ is  Gauduchon metric.  Choose  a $g$-Hermitian-Einstein metric $\lambda$ on ${\cal L}$ and denote by $C_\lambda$ the associated Chern connection. Therefore we have
  $$\frac{i}{2\pi} F_{C_\lambda}\wedge \omega_g^{n-1}=\frac{1}{\Vol_g}\deg_g({\cal L})vol_g\ ,
  $$
  where $vol_g=\frac{1}{n!}\omega_g^n$ is the volume form of $g$. Taking into account that $\omega_{FS}^k=0$ for $k\geq r$, this implies
  $$
  [\Omega_{g,h}]^{n+r-2}\wedge \frac{i}{2\pi} q^*\left[ F_{C_\lambda}\right]=\left(\begin{matrix}n+r-2\cr n-1\end{matrix}\right)  q^*(\frac{i}{2\pi}F_{C_\lambda}\wedge \omega_g^{n-1})\wedge [q^V_h]^*(\omega_{FS}^{r-1})=
  $$
  $$=\frac{1}{\Vol_g}\deg_g({\cal L}) (r-1)! \left(\begin{matrix}n+r-2\cr n-1\end{matrix}\right) q^*(vol_g)\wedge [q^V_h]^*(vol_{FS})=$$
 $$ \frac{1}{Vol_{H_{g,h}}(Z)}\frac{(n+r-2)!}{(n-1)!} \deg_g({\cal L})vol_{H_{g,h}}\ .
  $$
  This shows that the metric $q^*(\lambda)$ on the pull-back line bundle $q^*({\cal L})$ is Hermitian-Einstein with respect to $H_{g,h}$, and its Einstein constant is proportional to  $\deg_g({\cal L})$. But, in general one has (see Proposition 1.3.16 \cite{LT1}):

  \begin{re} 
 The Einstein constant of a Hermitian-Einstein metric on a holomorphic line bundle on a compact  Hermitian manifold $(M,H)$ is proportional  to the Gauduchon degree with respect to a Gauduchon metric in the conformal class of $H$.
 \end{re}
  
  Therefore, with respect to a Gauduchon metric $H'$ in the conformal class of $H_{g,h}$ one has
  \begin{equation}\label{volume}
vol_{H'}(D(\varphi))=  \deg_{H'}({\cal O}_Z(1)\otimes\pi^*({\cal L}^\vee))=\deg_{H'}({\cal O}_Z(1))- C\deg_g({\cal L})\ ,
 \end{equation}
for a positive constant $C$. Using this formula we get easily
\begin{pr}\label{compact} Let ${\cal E}$ be a holomorphic vector bundle on a compact complex  manifold $Y$ endowed with a Gauduchon metric $g$. For every $d\in\R$, the  subspaces 
$${\cal Q}uot^1_{\rm lf}({\cal E})_{\geq d}\subset {\cal Q}uot^1_{\rm lf}({\cal E})\ ,\ BN({\cal E})_{\geq d}\subset BN({\cal E})$$
 defined by the inequality $\deg_g({\cal L})\geq d$ are compact.
\end{pr}
\pf It suffices to prove the statement for  ${\cal Q}uot^1_{\rm lf}({\cal E})_{\geq d}$. But formula (\ref{volume}) shows that ${\cal Q}uot^1_{\rm lf}({\cal E})_{\geq d}$ can be identified with $[{\cal D}ou(Z)_{\leq s(d)}]\cap n_Z^{-1}(\im(a))$,  where  $s(d)$ has the form $s(d)=c_1-c_2 d$, for positive constants $c_1$, $c_2$. The result follows now from Proposition \ref{compactdouady}.
\qed
\begin{dt}
For a holomorphic bundle ${\cal E}$ on a Gauduchon manifold $(Y,g)$ we put
$$\degmax_g({\cal E}):=\left\{
\begin{array}{ccc} \sup\{\deg_g({\cal L})|\ [{\cal L}]\in BN({\cal E})\}&{\rm if}&BN({\cal E})\ne\emptyset\\
-\infty&{\rm if}&BN({\cal E})=\emptyset \ .
\end{array}\right.
$$
\end{dt}
\begin{re} \label{obtained} Let ${\cal E}$ be a holomorphic bundle over $Y$. 
 \begin{enumerate}
\item If $BN({\cal E})\ne\emptyset$ then the bound $\degmax_g({\cal E})$   is obtained, i.e. there exist  line bundles ${\cal L}$ on $Y$  with $H^0({\cal L}^\vee\otimes {\cal E})\ne 0$ and $\deg({\cal L})=\degmax_g({\cal E})$.
\item If $Y$ is projective  algebraic, then for every holomorphic bundle ${\cal E}$ on $Y$, one has  $BN({\cal E})\ne\emptyset$.
\item  For a rank 2 bundle ${\cal E}$ one has $BN({\cal E})\ne\emptyset$ if and only if ${\cal E}$ is a filtrable bundle.
\end{enumerate}
\end{re}
 \vspace{4mm}

The first statement follows from   Proposition \ref{compact}. Note that this statement is obvious when the pair $(Y,g)$ is algebraic (i.e.  when $g$ is a Hodge metric on $Y$). In this case $\deg_g$ takes values in $\Z$, so Remark \ref{obtained} (1) follows directly from  the boundedness  of the set $\deg_g(BN({\cal E}))$ (which is an obvious consequence of the formula (\ref{volume})). On the other hand, in the non-algebraic framework, and even in the K\"ahler case,  the image of the degree map $\deg_g:\Pic(Y)\to \R$ might be non-closed, so the boundedness of  the set $\deg_g(BN({\cal E}))$ is not sufficient.  
\\ 

Let now $X$ be an arbitrary complex manifold and let ${\cal Y}\to X$ be a complex manifold proper over $X$. Using Pourcin's relative Quot space \cite{Pou}, one can define the relative Douady space ${\cal D}ou_X({\cal Y})$, which is a complex space over $X$, whose fiber over a point $x\in X$ is ${\cal D}ou({\cal Y}_x)$. Proposition  \ref{compact} can be generalized to the relative case:
\begin{re}\label{compactpourcin} Let $h$ be a Hermitian metric on ${\cal Y}$. Then the topological subspaces
$${\cal D}ou_X({\cal Y})_{\leq s}:=\{D\in {\cal D}ou_X({\cal Y})|\ vol_h(D)\leq s\}\subset {\cal D}ou_X({\cal Y})
$$
are proper  over $X$.
\end{re}

Let ${\cal E}$ be a holomorphic bundle on  $X\times Y$ (a family of bundles on $Y$ parameterized by $X$).    Consider the {\it relative Brill-Noether locus } 
$$BN_X({\cal E}):=\{(x,[{\cal L}])|\ H^0({\cal L}^\vee\otimes{\cal E}^x)\ne 0\}\subset X\times \Pic(Y)\ ,
$$
which is obviously an analytic subset of $X\times \Pic(Y)$. We denote by ${\cal Q}uot^1_{{\rm lf},X}({\cal E})$ the open subspace of the relative Quot space ${\cal Q}uot_{X}({\cal E})$ (see \cite{Pou}) consisting of pairs $(x,q_x)$, where $q_x$ is a quotient of ${\cal E}_x$ with locally free rank 1 kernel. Set theoretically one has 
$${\cal Q}uot^1_{{\rm lf},X}({\cal E}):=\qmod{\{(x,\varphi)|\ x\in X,\ \varphi:{\cal L}\to {\cal E}^x\ne 0\}}{\sim}\ ,
$$
where the equivalence relation $\sim$ is induced by line bundle isomorphisms. The obvious projection on $X$ is  holomorphic. 
The same arguments as above and Remark \ref{compactpourcin}  give the following relative version of Proposition \ref{compact}:
\begin{pr}\label{rcompact} Let ${\cal E}$ be a holomorphic vector bundle on the product $X\times Y$ and $g$ a Gauduchon metric on  $Y$. For every $d\in\R$, the  subspaces  
$${\cal Q}uot^1_{{\rm lf},X}({\cal E})_{\geq d}\subset {\cal Q}uot^1_{{\rm lf},X}({\cal E})\ ,\ BN_X({\cal E})_{\geq d}\subset BN_X({\cal E})$$
 defined by the inequality $\deg_g({\cal L})\geq d$ are proper over $X$, so their projections  $X_{\geq d}$ on $X$  are closed with respect to the classical topology.
\end{pr}

In particular ${\cal Q}uot^1_{{\rm lf},X}({\cal E})_{\geq d}$, $BN_X({\cal E})_{\geq d}$ are compact when $X$ is compact. In this case the degree map $\deg_g$ will be  bounded on  any $BN_X({\cal E})_{\geq d}$, so  
\begin{re}\label{empty} If $X$ is compact, then $ BN_X({\cal E})_{\geq d}$ and ${\cal Q}uot^1_{{\rm lf},X}({\cal E})_{\geq d}$  are empty for sufficiently large $d$.
\end{re}
In general, the subsets $X_{>d}$ defined in a similar way are not necessary closed even with respect to the classical topology.  Note that $X_{\geq d}$, $X_{> d}$ can be defined as
$$X_{\geq d}=\{x\in X|\ \degmax({\cal E}_x)\geq d\}\ ,\ X_{>d}=\{x\in X|\ \degmax({\cal E}_x) > d\}\ .
$$
\begin{re}\label{sets} It holds
\begin{enumerate}
\item $$X_{> d'}\subset X_{\geq d'}\subset X_{> d}\subset X_{\geq d}\hbox{ for } d' > d\ ;\ \bigcap_{d \geq a} X_{\geq d}=\emptyset\ .$$ 
\item When $\rk({\cal E})=2$, then $\union_{d\leq a} X_{\geq d}={\rm Filtr}({\cal E})$,
where ${\rm Filtr}({\cal E})$ denotes the set of points $x\in X$ for which ${\cal E}_x$ is filtrable. This set is not closed in general (counterexamples are known when $Y$ is a class VII surface with $b_2=1$).
\item When $X$ is compact, then the set $\{d\in\R|\ X_d\ne\emptyset\}$ is bounded from above and, when non-empty, has a maximum, which coincides with 
$$\max\{\deg_g({\cal L})|\ (x,{\cal L})\in BN_X({\cal E})\}= \max\{\degmax_g({\cal E}_x)|\ \ x\in X\}\ .$$
\end{enumerate}
\end{re}

\vspace{4mm}

 Denote by $p_\Pic$ the projection $X\times\Pic(Y)\to \Pic(Y)$. Proposition \ref{rcompact} has important consequences in the particular case when $\deg_g\circ p_\Pic$ is locally constant on   $BN_X({\cal E})$. This happens automatically when $\deg_g$ itself is locally constant (i.e. it is a topological invariant).
\begin{pr} \label{topinv} Suppose that $\deg_g$ is locally constant on $BN_X({\cal E})$.  Then    %
\begin{enumerate}
\item ${\cal Q}uot^1_{{\rm lf},X}({\cal E})_{\geq d}$ is closed and open in  ${\cal Q}uot^1_{{\rm lf},X}({\cal E})$, whereas  $BN_X({\cal E})_{\geq d}$ is  closed and open in $BN_X({\cal E})$. 
\item ${\cal Q}uot^1_{{\rm lf},X}({\cal E})_{\geq d}$,  $BN_X({\cal E})_{\geq d}$ are complex spaces which are proper over $X$, 
\item The open subspace ${\cal Q}uot^1_{{\rm lf},X}({\cal E})_{>d}$ ($BN_X({\cal E})_{>d}$) is closed in ${\cal Q}uot^1_{{\rm lf},X}({\cal E})_{\geq d}$ (respectively $BN_X({\cal E})_{\geq d}$), so it is   proper over $X$, too.
\item The  projections $X_{\geq d}$, $X_{> d}$ of these spaces on $X$ are Zariski closed.
\end{enumerate}
\end{pr}
Proposition \ref{topinv} applies for instance to families parameterized by a compact manifold:

\begin{re}\label{compactbase} Suppose that the parameter manifold $X$  is compact. Then $\deg_g\circ p_\Pic$ is locally constant on $BN_X({\cal E})$.
\end{re}
\pf Let $C$ be a connected component of $BN_X({\cal E})$. Choose a point $z_0=(x_0,{\cal L}_0)\in C$ and $d<\deg_g({\cal L}_0)$. The subset $C_{\geq d}:=C\cap BN_X({\cal E})_{\geq d}$ is non-empty and compact, so  $\deg_g\circ p_\Pic$ reaches its maximum on this set. This maximum, say $m$, is also the maximum of $\resto{\deg_g\circ p_\Pic}{C}$. Note that the degree map $\deg_g:\Pic(Y)\to \R$ is pluriharmonic (indeed,  its restriction to $\Pic^0(Y)$ is induced by a  an $\R$-linear map $H^1(Y,{\cal O}_Y)\to\R$). Therefore the nonempty level set
$$C_{=m}:=\{z\in C|\ \deg_g\circ p_\Pic(z)=m\}
$$ 
(which is obviously closed) is also open in $C$. To see this it suffices to apply the usual open mapping theorem  to the induced pluriharmonic function on a desingularization of $C$. Therefore $C_{=m}=C$, which   completes the proof.\qed
\begin{re} \label{oldremark} Suppose that   $X$ is compact. Then  the set of jumps 
$$J:=\{\degmax({\cal E}_x)|\ x\in X\}\setminus\{-\infty\}\subset\R$$
 is either finite or it forms a strictly decreasing sequence which tends to $-\infty$.   
When the family has the property that any bundle ${\cal E}_x$ in the family admits invertible subsheaves (in particular when any ${\cal E}_x$ is filtrable) then $J$ is finite.
\end{re}
\pf  We know by Remark \ref{compactbase} that $\deg_g\circ p_\Pic$ is locally constant on $BN_X({\cal E})$, so $BN_X({\cal E})_{\geq d}$ is a compact complex space for any  $d\in\R$. Note that
$$J\cap [d,\infty)\subset \deg_g\circ p_\Pic(BN_X({\cal E})_{\geq d})\ .
$$
But $\deg_g\circ p_\Pic$ takes only finitely many values on the compact complex space $BN_X({\cal E})_{\geq d}$. Therefore $J\cap [d,\infty)$ is finite for every $d\in\R$, which proves the first claim.

 For the second, note that the hypothesis implies  $\union_{d\leq a} X_{d}=X$ (see Remark \ref{sets}). It suffices to note that a compact complex space (or more generally a complex space with finitely many irreducible components) cannot be written as the a union of  an infinite countable chain of analytic subsets with strict inclusions.
\qed

 \subsection{Openness properties}
 
 Let $V$ be a complex vector space of dimension $r$ and let $s\in\{1,\dots,r-1\}$. The Grassmannian $G_s(V)$ can be naturally embedded in $\P(\wedge^s(V))$ using the Pl\"ucker embedding 
 $$\langle v_1,\dots, v_s\rangle\mapsto [v_1\wedge\dots\wedge v_s]\ .$$
 The cone $C_s(V)\subset \wedge^s(V)$ over the Grassmannian $G_s(V)$ will be called the cone of exterior monomials in $\wedge^s(V)$.   For an element  $u\in\wedge^s(V)$ we denote by $u^\bot$ the subspace of $V^*$ consisting of forms $v$ with $\iota_v(u)=0$ and by $(u^\bot)^\bot\subset V$ the annihilator of $u^\bot$. Note that the map 
 $$C_s(V)\setminus\{0\}\ni u\mapsto (u^\bot)^\bot\subset V$$
 is just the standard projection on the Grassmannian $G_s(V)=\P(C_s(V))$.

 Similarly, for a holomorphic rank $r$ bundle ${\cal E}$ over $Y$  we obtain a locally trivial, closed cone sub-bundle $C_s({\cal E})\subset \wedge^s({\cal E})$ over $Y$. If ${\cal L}\subset\wedge^s{\cal E}$ is a rank 1 locally free subsheaf which is contained in $C_s({\cal E})$, then  $[{\cal L}^\bot]^\bot$ will be a rank $s$ subsheaf of ${\cal E}$ which is a subbundle on the open set where ${\cal L}$ is a subbundle.
 
 The following result is inspired by L\"ubke's proof of the "first implication"  of the Kobayashi-Hitchin correspondence (see \cite{Lu}, \cite{LT1}). This result is useful because it shows that (semi)stability can be tested using only rank 1 subsheaves of a finite set of bundles associated with ${\cal E}$.
 
\begin{pr}\label{martin} The following conditions are equivalent:
\begin{enumerate}
\item ${\cal E}$ is $g$-stable ($g$-semistable).
\item For every $s\in\{1,\dots,r-1\}$ and any non-trivial morphism $\varphi:{\cal L}\to \wedge^s({\cal E})$ defined on a line bundle ${\cal L}$ and having the property $\im(\varphi)\subset C_s({\cal E})$ one has
$$\deg_g({\cal L})< s\mu_g({\cal E})\ \hbox{ (respectively } \deg_g({\cal L})\leq s\mu_g({\cal E}) )\ .
$$
\end{enumerate}
\end{pr}
\pf We treat only the equivalence for stability.
\\ 
\\
(2) $\Rightarrow$ (1): Let  ${\cal F}\subset {\cal E}$ be a rank $s$ subsheaf   with torsion-free quotient. We obtain a morphism
$$\det({\cal F})=\Lambda^s({\cal F})^{**}\to \wedge^s({\cal E})^{**}=\wedge^s({\cal E})\ ,
$$
which is a bundle embedding on the complement of  a codimension 2 analytic subset, and whose image is contained in $C_s({\cal E})$. By (2)  one must have 
$$s\mu_g({\cal F})=\deg_g(\det({\cal F}))<s\mu_g({\cal E})\ ,$$
which shows that the stability condition is satisfied.
\\
\\
(1)$\Rightarrow$ (2): Let $\varphi:{\cal L}\to \wedge^s({\cal E})$ be a non-trivial morphism whose image is contained in $C_s({\cal E})$. Consider the exact sequence
$$0\map {\cal L}\stackrel{\varphi}{\hookrightarrow}  \wedge^s({\cal E})\textmap{q} {\cal Q}\map 0
$$
associated with $\varphi$. The sheaf  ${\cal L}':=q^{-1}(\Tors(Q))$ is a rank 1 reflexive sheaf, hence it is a line bundle. Since $\varphi$ induces a sheaf monomorphism  ${\cal L}\hookrightarrow  {\cal L}'$, we get   $\deg_g({\cal L}')\geq\deg_g({\cal L})$. Note that ${\cal L}'$ is still contained in $C_s({\cal E})$ and the sheaf inclusion ${\cal L}'\subset \wedge ^s({\cal E})$ is a bundle embedding on the complement $U$ of  a codimension 2 analytic subset.  Therefore, the subsheaf ${\cal F}:=([{\cal L}']^\bot)^\bot$ of ${\cal E}$ is a rank $s$ {\it subbundle}  on $U$. The line bundles ${\cal L}'$ and $\det({\cal F}')$ are isomorphic (by Hartogs theorem, because they are isomorphic on $U$). Therefore, using the hypothesis (${\cal E}$ is stable), one gets
$$\deg({\cal L})\leq \deg_g({\cal L}')=\deg_g({\cal F})< s\mu_g({\cal E})\ .$$
\qed

We can prove now the openness properties of the stability condition which have been announced in the introduction (Theorem \ref{open}): \\
\\
 \pf (of Theorem \ref{open})
Suppose first that $r=2$. In this case
  $$X^\st=X\setminus X_{\geq \mu_g({\cal E})}\ ,\  X^\sst=X\setminus X_{>\mu_g({\cal E})}\ .
 $$
 Therefore $X^\st$ is open by Proposition \ref{rcompact} for general families, and $X^\st$, $X^\sst$ are Zariski open by Remark \ref{topinv} as soon as the map $\deg_g\circ p_\Pic$ is locally constant on $BN_X({\cal E})$. This happens when $\deg_g$ is a topological invariant or when $X$ is compact (see Remark \ref{compactbase}). For the case $r>2$, we make use of Proposition \ref{martin}. Let 
 $$BN_X^s({\cal E}):=\{(x,[{\cal L}])|\ \exists \varphi:{\cal L}\to \wedge^s({\cal E}_x) \hbox{ with } \varphi\ne 0,\ \im(\varphi)\subset C_\varphi({\cal E}) \}\subset BN_X (\wedge^s({\cal E})\ ,
$$
and let  ${\cal Q}uot^1_{{\rm lf},X}({\cal E})^s$ be the closed complex subspace of  ${\cal Q}uot^1_{{\rm lf},X}(\wedge^s({\cal E}))$ consisting of quotients with kernel contained in $C_s({\cal E})$.  We obtain subsets $X^s_{\geq d}$ which are closed for general families, and subsets $X^s_{\geq d}$, $X^s_{> d}$ which are Zariski closed when one of the assumptions (\ref{ti}), (\ref{cb}) holds (because in these cases $\deg_g\circ p_\Pic$ will be locally constant on $BN_X^s({\cal E})$).
On the other hand, by Proposition \ref{martin} we get
$$X^\st=X\setminus \union_{1\leq s\leq r-1} X^s_{\geq s\mu_g({\cal E})} \ , X^\sst=X\setminus \union_{1\leq s\leq r-1} X^s_{>s\mu_g({\cal E})}
$$
which completes the proof.
\qed

\begin{re}
\begin{enumerate}
\item  The first statement in Proposition \ref{open} can be proved using the Kobayashi-Hitchin correspondence and the implicit function theorem \cite{LT1}.  
\item For general families $X^\st$ is not necessary Zariski open and $X^\sst$ is not necessary open even in the classical topology (see the example below). 
\end{enumerate}
\end{re}
\noindent
{\bf Example:} (see \cite{Te1}, \cite{Te2} for details) Let $Y$ be a minimal class VII surface with $b_2=1$ and no homologically trivial divisor. For every ${\cal L}\in\Pic^0(Y)$ we have an (essentially unique)   non-trivial extension 
$$0\map {\cal L}\map {\cal E}_{\cal L}\map {\cal K}\otimes {\cal L}^{-1}\map 0\ .
$$
The bundles ${\cal E}_{\cal L}$ form a holomorphic family parameterized by $X=\Pic^0(Y)\simeq\C^*$.   ${\cal E}_{\cal L}$ is stable if and only if $\deg_g({\cal L})<\frac{1}{2}\deg({\cal K})$. This inequality defines a pierced open disk, which is obviously open but not Zariski open in the parameter space $X$.
The set $X^\sst$ is the pierced closed disk defined by the non-strict inequality $\deg_g({\cal L})\leq \frac{1}{2}\deg({\cal K})$, so it is not open in the classical topology.

\parindent=12 pt
\section{The Petersson-Weil current on the base of a  generically stable family}
\label{PW}

\subsection{Donaldson's functional for families. An integral formula}

Let $(Y,g)$ be a compact $n$-dimensional K\"ahler manifold and $\omega$ its K\"ahler form. Let ${\cal E}$ be a rank $r$ holomorphic bundle on $Y$ with determinant ${\cal L}:=\det({\cal E})$, and let $l$ be a fixed Hermitian metric on ${\cal L}$.  We consider the space ${\cal M}et^l({\cal E})$ of Hermitian metrics $H$ on ${\cal E}$ with $\det(H)=l$. The bundle ${\cal E}$ is polystable if and only it admits a projectively Hermitian-Einstein metric $H\in {\cal M}et^l({\cal E})$, i. e. a   metric $H\in {\cal M}et^l({\cal E})$ which  satisfies the projective Hermitian-Einstein equation:
$$i\Lambda F_H^0=0\ ,
$$
where $F_H^0$ denotes the trace-free part of the curvature $F_H$ of the Chern connection $A_H$ of $H$. This equation is equivalent with to  the weak Hermitian-Einstein equation 
$$i\Lambda F_H=c\id_{\cal E}\ .
$$
with Einstein factor $c:=\frac{1}{r}i\Lambda F_l$.  If ${\cal E}$ is stable then it admits a {\it unique} projectively Hermitian-Einstein metric $H_{he}\in {\cal M}et^l({\cal E})$.

Let 
$$M(\cdot,\cdot):{\cal M}et^l({\cal E})\times {\cal M}et^l({\cal E})\map \R$$
 be the (projective) Donaldson functional on the space of pairs of metrics. For two metrics $K$, $H\in{\cal M}et^l({\cal E})$  the real number $M(K,H)$ is defined by
 \begin{equation}\label{Mdef}
 M(K,H)=\int_Y R_2^0(K,H)\wedge \omega^{n-1}\ ,
\end{equation}
 where 
 $$R_2^0(K,H)\in \qmod{A^{1,1}_\R(Y)}{(\partial A^{0,1}(Y)\oplus \bar\partial A^{1,0}(Y))\cap A^{1,1}_\R(Y)}
 $$
 is the Bott-Chern secondary holomorphic characteristic class   of the pair $(K,H)$ associated with the ad-invariant symmetric function 
 $$ (a,b)\mapsto \tr(a^0b^0)$$
  on the Lie algebra $gl(r,\C)$. One has (see \cite{Do1})
 \begin{equation}\label{ddbar}
 i\bar\partial\partial R_2^0(K,H)=\tr (F_K^0\wedge F_K^0)-\tr (F_H^0\wedge F_H^0)\ .
 \end{equation}
  The assignments    $R_2^0(\cdot,\cdot)$, $M(\cdot,\cdot)$  satisfy the following variational formulae:
\begin{equation}\label{variation}
\frac{d}{dt} R_2^0(K,H_t)=2i\tr(H^{-1}_t\dot{H}_t F_{H_t}^0)\ ,$$
$$\frac{d}{dt} M(K,H_t)=2\int\limits_Y i\tr(H^{-1}_t\dot{H}_t F_{H_t}^0)\wedge \omega^{n-1}\ .
 \end{equation}

Let now $X$ be a  connected complex manifold (not necessary compact) and ${\cal E}$ a holomorphic bundle over the product $X\times Y$, i.e. a family of holomorphic bundles on $Y$ parameterized by $X$.    Let $l$ be a fixed Hermitian metric on ${\cal L}:=\det ({\cal E})$. For every Hermitian metric $H\in {\cal M}et^l({\cal E})$, we consider the Chern form
$$\eta_H:=\tr (F_H^0\wedge  F_H^0)=\frac{4\pi^2}{r}\left[2rc_2({\cal E},H)-(r-1)c_1({\cal E},H)^2\right]\ ,$$
where $c_i({\cal E},H)$ denotes the $i$-th Chern form of the Chern connection associated with the pair $({\cal E},H)$. Let $H$, $K\in {\cal M}et^l({\cal E})$. We introduce the function  $m(K,H)\in{\cal C}^\infty(X,\R)$ defined by
$$m(K,H)(x)=M(K_x,H_x)\ ,
$$
where $M_Y(\cdot,\cdot)$ denotes the Donaldson functional on the space  ${\cal M}et^{l_x}({\cal E}_x)\times  {\cal M}et^{l_x}({\cal E}_x)$ of pairs of metrics with fixed determinant $l_x$ on  the bundle ${\cal E}_x$ over $Y$. The function $m(K,H)$ will be called the {\it fiberwise Donaldson functional} of the pair $(K,H)$.

\begin{pr} \label{integral} Let $K$, $H\in {\cal M}et^l({\cal E})$  be two Hermitian metrics on ${\cal E}$. Then
$$[p_X]_*\left[\eta_H\wedge p_Y^*(\omega_Y^{n-1})- \eta_K\wedge p_Y^*(\omega_Y^{n-1})\right]=i\partial\bar\partial m(K,H)
$$
\end{pr}
\pf Using formula (\ref{ddbar}) one obtains
$$[p_X]_*\left[\eta_H\wedge p_Y^*(\omega_Y^{n-1})- \eta_K\wedge p_Y^*(\omega_Y^{n-1})\right]=[p_X]_*(i\partial\bar\partial R_2^0(K,H)\wedge\omega_Y^{n-1})=$$
$$=i\partial\bar\partial\left[[p_X]_* (R_2^0(K,H)\wedge \omega_Y^{n-1})\right]\ .
$$
It suffices to use (\ref{Mdef}) and the obvious functoriality property of the Bott-Chern secondary holomorphic characteristic classes with respect to holomorphic maps, which implies in our case:
$$\resto{R_2^0(K,H)}{\{x\}\times Y}=R_2^0(K_x,H_x)\ .
$$
\qed
 
The background  metric $K$ should be considered as a known object. In particular  suppose that a local potential $\varphi_K\in{\cal C}^\infty(U,\R)$  ($U\subset X$) for the closed (1,1)-form $ [p_X]_*(\eta_K\wedge p_Y^*(\omega_Y^{n-1}))$ is known. Proposition \ref{integral} furnishes a local potential of  the closed (1,1)-form $ [p_X]_*(\eta_H\wedge p_Y^*(\omega_Y^{n-1}))$ associated with the variable metric $H$ in terms of the known local potential  $\varphi_K$ and the fiberwise Donaldson functional $m(K,H)$:
\begin{re}\label{mainremark} The  function $\varphi_K+ \resto{m(K,H)}{U}$ is a local potential for the closed (1,1)-form $ [p_X]_*(\eta_H\wedge p_Y^*(\omega_Y^{n-1}))$.
\end{re}
A holomorphic rank $r$-vector bundle ${\cal E}$ defines a holomorphic principal $GL(r,\C)$-bundle ${\cal P}_{\cal E}$ (the frame bundle of ${\cal E}$) and a holomorphic $PGL(r,\C)$-bundle 
$${\cal Q}_{\cal E}:=\qmod{{\cal P}_{\cal E}}{\C^*}\ .$$

The stability conditions for ${\cal E}$, ${\cal P}_{\cal E}$ and ${\cal Q}_{\cal E}$ are equivalent. Consider the open  subset  
$$X^{\rm st}:=\{x\in X|\ {\cal E}_x\hbox{ is stable}\}\subset X\ ,
$$
which is Zariski open, by Theorem \ref{open}. Let $E$ be a differentiable vector bundle which is isomorphic to the underlying differentiable bundles of the bundles ${\cal E}_x$ and let $Q:=P_E/\C^*$ be the associated differentiable principal $PGL(r,\C)$-bundle. Denote by
${\cal M}^\st(Q)$   the moduli space of stable holomorphic structures on $Q$. In the definition of this moduli space it is convenient to use the ``small" gauge group $\Aut_0(Q):=\Gamma(Y,Q\times_{\rm Ad}SL(r,\C))$ instead  of $\Aut(Q)= \Gamma(Y,Q\times_{\rm Ad}PGL(r,\C))$, because in this way one obtains a moduli space which is naturally isomorphic with the moduli space ${\cal M}^\st_{\cal D}(E)$ of stable holomorphic structures on $E$ which induce a fixed holomorphic structure ${\cal D}$ on $\det(E)$ (which therefore is independent of ${\cal D}$ up to canonical isomorphism). The image of $\Aut_0(Q)$ in $\Aut(Q)$ has finite index \cite{LT1}, \cite{Te1}. ${\cal M}^\st(Q)$ is a (possibly singular) K\"ahler space. The corresponding Hermitian structure on the Zarisky tangent space $T_{\cal Q}({\cal M}^\st(Q))=H^1(\ad({\cal Q}))$ is given by the $L^2$-product on the harmonic space $\H^{0,1}_{B}(\ad({\cal Q}))$ with respect to the Hermite-Einstein connection $B$ on $Q$ associated with the unique Hermite-Einstein $PU(r)$-reduction of the holomorphic structure ${\cal Q}$  \cite{LT2}.

Our family $\resto{{\cal E}}{X^\st\times Y}$ induces a holomorphic map $F:X^{\rm st}\to {\cal M}^\st$  given by
$$F(x)=[{\cal Q}_{{\cal E}_x}]\ .
$$
 Solving fiberwise the projective Hermitian-Einstein equation, one gets a {\it smooth} Hermitian metric $H_{he}$ on $X^{\rm st}\times Y$.

 \begin{pr} \label{PWform} Let $\Omega$ be the K\"ahler form of the Petersson-Weil metric on the K\"ahler space ${\cal M}^\st(Q)$. Then  
\begin{equation}\label{ST}
2F^*(\Omega)=[p_{X^\st}]_*(\eta_{H_{he}}\wedge p_Y^*(\omega_Y^{n-1})) \ .
 \end{equation}
 \end{pr}
\pf The bundle $\Lambda^1_{X\times Y}$ splits as $\Lambda^1_X\oplus\Lambda^1_Y$, where  $\Lambda^1_X$, $\Lambda^1_Y$ denote the pull-backs of the corresponding bundle of forms on $X$ and $Y$ via the two projections. Similarly the bundle $\Lambda^1_{X\times Y}\otimes E$ of $E$-valued 1-forms decomposes as $\Lambda^1_{X\times Y}\otimes E=\Lambda^1_X\otimes  E\oplus  \Lambda^1_Y\otimes  E$. Taking into account the complex structures, we can decompose further
\begin{equation}\label{decforms}
\Lambda^1_X\otimes E=\Lambda^{X,0}(E)\oplus \Lambda^{0,X}(E)\ ,\  \Lambda^1_Y\otimes E=\Lambda^{Y,0}(E)\oplus \Lambda^{0,Y}(E)\ ,\ 
\end{equation}
where
$$\Lambda^{X,0}(E)=\Lambda^{1,0}_X\otimes E,\ \Lambda^{0,X}(E)=\Lambda^{0,1}_X\otimes E,\ \Lambda^{Y,0}(E)=\Lambda^{1,0}_Y\otimes E,\ \Lambda^{0,Y}(E)=\Lambda^{0,1}_Y\otimes E.
$$
The bundles of $E$-valued $p$-forms split similarly.  The covariant derivative $\nabla_A$  of a connection  $A$ on $E$  decomposes as 
$$\nabla_{A}^X=\nabla_{A}^{X,0}+\nabla_{A}^{0,X}\ ,\nabla_{A}^Y=\nabla_{A}^{Y,0}+\nabla_{A}^{0,Y}
$$
We denote by $d_{A}^{X,0}$, $d_{A}^{0,X}$, $d_{A}^{Y,0}$, $d_{A}^{0,Y}$ the  natural extensions of these operators on the spaces of forms  $A^{*}(E)$,  $A^{*}(\End(E))$. Suppose now that $A=A_H$ is the integrable Hermitian connection associated with a Hermitian metric $H$ on ${\cal E}$. The curvature $F_H:=F_{A_H}$ has type $(1,1)$, hence it decomposes as 
$$F_H=F_H^{X,X}+F_H^{Y,Y}+F_H^{X,Y}+ F_H^{Y,X}\ ,
$$
where 
$$F_H^{X,X}=d_{A_H}^{X,0}\circ d_{A_H}^{0,X}+d_{A_H}^{0,X}\circ d_{A_H}^{X,0}\ ,\ F_H^{Y,Y}=d_{A_H}^{Y,0}\circ d_{A_H}^{0,Y}+d_{A_H}^{0,Y}\circ d_{A_H}^{Y,0}\ , $$
$$F_H^{X,Y}=d_{A_H}^{X,0}\circ d_{A_H}^{0,Y}+d_{A_H}^{0,Y}\circ d_{A_H}^{X,0}\ ,\ F_H^{Y,X}=d_{A_H}^{Y,0}\circ d_{A_H}^{0,X}+d_{A_H}^{0,X}\circ d_{A_H}^{Y,0}\ .
$$
whereas the (2,0)+(0,2)   components 
$$ F_{H}^{XY,0}=d_{A_H}^{X,0}\circ d_{A_H}^{Y,0}+d_{A_H}^{Y,0}\circ d_{A_H}^{X,0}\ ,\ F_{H}^{0,XY}=d_{A_H}^{0,X}\circ d_{A_H}^{0,Y}+d_{A_H}^{0,Y}\circ d_{A_H}^{0,X}\ , $$
$$F_{H}^{XX,0}=d_{A_H}^{X,0}\circ d_{A_H}^{X,0}\ ,\  F_{H}^{YY,0}=d_{A_H}^{Y,0}\circ d_{A_H}^{Y,0}\ , \ F_{H}^{0,XX}=d_{A_H}^{0,X}\circ d_{A_H}^{0,X}\ ,\ F_{H}^{0,YY}=d_{A_H}^{0,Y}\circ d_{A_H}^{0,Y}
$$
vanish. Taking into account  these relations, the Bianchi identity $d_A F_A=0$ shows that all the terms
\begin{equation}\label{bi}
d_{A_H}^{X,0}(F_H^{X,Y}),\ d_{A_H}^{0,Y}(F_H^{X,Y}),\ d_{A_H}^{Y,0}(F_H^{Y,X}),\ d_{A_H}^{0,X}(F_H^{Y,X}),\ d_{A_H}^{X,0}(F_H^{Y,X})+d_{A_H}^{Y,0}(F_H^{X,X}),$$
$$ d_{A_H}^{0,Y}(F_H^{Y,X})+d_{A_H}^{0,X}(F_H^{Y,Y}),\ d_{A_H}^{Y,0}(F_H^{X,Y})+d_{A_H}^{X,0}(F_H^{Y,Y}),\ d_{A_H}^{0,X}(F_H^{X,Y})+d_{A_H}^{0,Y}(F_A^{X,X})
\end{equation}
must vanish, too. For every  tangent vector $v\in T^{1,0}_x(X)$ one obtains  a $\End(E_x)$-valued type $(0,1)$-form $h(v):=\iota_v  F_H^{X,Y}$ on $\{x\}\times Y\simeq Y$, and the  vanishing of $d_{A_H}^{0,Y}(F_A^{X,Y})=0$ shows that $h(v)$ is $\bar\partial_{{\cal E}_x}$-closed.  If one uses the connection $A_{H'}=A+\alpha^{1,0}$ associated with a new metric $H'$, the corresponding form $h'(v)$ will be $h(v)+\bar\partial_{{\cal E}_x}(\iota_v(\alpha)))$, so it will define the same Dolbeault cohomology class $\chi(v)$ as $h(v)$. It is well known that $\chi(v)$ is precisely the infinitesimal deformation  of the holomorphic bundle ${\cal E}_x$ in the direction $v$ associated with the family ${\cal E}$ (see \cite{ST}  Proposition 1 p.   102).  The corresponding infinitesimal deformation $F_*(v)$ of the $PGL(r,\C)$-bundle ${\cal Q}_{{\cal E}_x}$  will be $\chi_0(v):=[h_0(v)]$, where $h_0(v):=\iota_v  [F_H^{X,Y}]^{0}$. The vanishing of the seventh term in (\ref{bi}) shows that
$$\Lambda_Y \partial_{A_{H_x}} h(v)=\nabla_{A_H,v}\Lambda_Y  F^{Y,Y}_{H_x} \ ,\ \Lambda_Y \partial_{A_{H_x}} h_0(v)=\nabla_{A_H,v}\Lambda_Y [F^{Y,Y}_{H_x}]^0\ .$$
So $\Lambda_Y \partial_{A_{H_x}} h_0(v)$ vanishes when $H$ is fiberwise projectively Hermite-Einstein around the fiber $\{x\}\times Y$. Therefore, in this case, $h_0(v)$ is just the harmonic  representative of $F_*(v)\in H^1(\ad({\cal Q}_x))$ with respect to the Hermite-Einstein $PU(r)$-reduction of ${\cal Q}_x$ defined by $H_x$. 
One has
$$[p_{X}]_*(\eta_{H}\wedge p_Y^*(\omega_Y^{n-1})) =$$
$$=2[p_{X}]_*\left\{ \tr \left[ [F^{X,Y}_{H}]^0\wedge  [F^{Y,X}_{H}]^0 +[F^{X,X}_{H}]^0\wedge  [F^{Y,Y}_{H}]^0\right]\wedge p_Y^*(\omega_Y^{n-1})\right\}$$
If  now $H$ is fiberwise projectively Hermite-Einstein, the second term on the right vanishes, and for an vector $v\in T^{1,0}_X$ one obtains
$$[p_{X}]_*\left\{ \tr \left[ [F^{X,Y}_{H}]^0\wedge  [F^{Y,X}_{H}]^0 \right]\wedge p_Y^*(\omega_Y^{n-1})\right\}(v,\bar v)=i \| h_0(v)\|^2=\Omega(F_*(v),F_*(\bar v))\ .
$$
 \qed

According to \cite{ST} the closed positive (1,1)-form $[p_{X^\st}]_*(\eta_{H_{he}}\wedge p_Y^*(\omega_Y^{n-1}))$ on $X^\st$ will be called the (projective) {\it Petersson-Weil form} of the family. We will use the notation $pw_0(\resto{{\cal E}}{X^\st\times Y})$ for this  form.

\begin{pr} \label{neg} Suppose that ${\cal E}$ is a stable bundle on $(Y,g)$ and let $H_{he}\in {\cal M}et^l({\cal E})$ be its unique projective Hermitian-Einstein metric. Then for any $K\in {\cal M}et^l({\cal E})$ one has
\begin{equation}\label{negeq}
M(K,H_{he})\leq 0 .
\end{equation}
Equality occurs if only   if $K=H_{he}$.
\end{pr}
\pf \footnote{The simple argument given here has been kindly suggested by the referee.} Let $S$  be a trace-free $K$-Hermitian endomorphism of ${\cal E}$, and put $H_t:=K e^{tS}$.  Using the second variational formula in (\ref{variation}) one gets
$$\frac{d}{dt} M(K,H_t)=2\int\limits_Y i\tr(S F_{H_t}^0)\wedge \omega^{n-1}\ ,\ $$
$$\frac{d^2}{dt^2} M(K,H_t)=2\int\limits_Y i\tr(S \bar\partial_{\cal E}\partial_{H_t} S)\wedge \omega^{n-1}=2\| \partial_{H_t} S\|_{L^2}^2=2\| \bar\partial_{\cal E} S\|_{L^2}^2\ .
$$
Therefore, the function $f_{K,S}:t\mapsto M(K,K e^{tS})\in\R$ is convex. Choose now $S$ such that   $H_{he}=K e^{S}$. Using the first formula above and the Hermite-Einstein equation, we get 
$$\at{\frac{d}{dt}}{t=1}  M(K,H_t)=0\ ,$$
so the convex function $f_{K,S}$ is non-increasing on $[0,1]$.  But $f_{K,S}(0)=0$, hence $f_{K,S}(1)=M(K,H_{he})\leq 0$. The second statement follows easily using the fact that stable bundles are simple.
\qed

We now come back to a holomorphic bundle ${\cal E}$ on   $X\times Y$ (where $X$ is a connected complex manifold) endowed with a fixed Hermitian metric $l$ on its determinant line bundle, and a background Hermitian metric $K\in {\cal M}et^l({\cal E})$.  Let $H_{he}\in {\cal M}et^l(\resto{{\cal E}}{X^\st\times Y})$ be the Hermitian metric obtained by solving fiberwise the projective Hermite-Einstein equation. 
 \begin{co} \label{smallerpotential} Let $\varphi_K\in{\cal C}^\infty(U,\R)$ be a local potential of   $[p_X]_*(\eta_K\wedge p_Y^*(\omega_Y^{n-1}))$ defined on $U\subset X$. Then   the local potential   
 \begin{equation}
 \label{smallerpotentialPW}
 \psi:U\cap X^\st\to \R\ ,\ \psi:=\resto{\varphi_K}{U\cap X^\st}+\resto{m(K,H_{he})}{U\cap X^\st}
 \end{equation}
of the Petersson-Weil form $pw_0(\resto{{\cal E}}{X^\st\times Y})$ satisfies the inequality
 \begin{equation}\label{smallerPW}
 \psi\leq \resto{\varphi_K}{U\cap X^\st}\ .
 \end{equation}
  \end{co} 
  \pf This follows directly from Remark \ref{mainremark} and Proposition  \ref{PWform}.
  \qed
  
  We can prove now the extensibility result  stated in the introduction: When $X^\st\ne\emptyset$, the Petersson-Weil form on $X^\st$ extends as a closed positive current on $X$.
  \\
  \\
  \pf (of Theorem \ref{extension}) In general, a plurisubharmonic function defined on the complement of an analytic set $E$ in a complex manifold $Z$ extends as a plurisubharmonic function on $Z$ if it is locally bounded {\it from above} near $E$ (see for instance Theorem 6.2.9 in \cite{DS}).  Therefore, by (\ref{smallerPW}),  a local potential $\psi:U\cap X^\st\to \R$ obtained using Corollary \ref{smallerpotential} extends as a plurisubharmonic function $\tilde\psi$  on $U$ (recall that $X^\st$ is  {\it Zariski open} by Theorem \ref{open}, so its complement is an analytic set). 
  Two potentials $\psi$, $\psi'$ associated with  potentials $\varphi_K:U\to\R$, $\varphi'_K:U'\to\R$ of $[p_X]_*(\eta_K\wedge p_Y^*(\omega_Y^{n-1}))$ differ by a  pluriharmonic function which is defined on $U\cap U'$ (this follows from (\ref{smallerpotentialPW})). Therefore $\tilde\psi$, $\tilde\psi'$ define the same positive current on $U\cap U'$. In this way we obtain a global current on $X$, which does not depend on  the background metric $K$. 
\qed
\begin{dt}  The closed positive current given by Theorem \ref{extension} will be denoted  by $pw_0({\cal E})$ and will be called  the {\it  projective Petersson-Weil current of the family ${\cal E}$.}
\end{dt}
\begin{re}\label{interesting} \begin{enumerate} 
\item A  positive current defined on the complement of an analytic set $A$ may have many extensions, because one can add Dirac   type currents  concentrated on the irreducible components of $A$. Our result shows that $pw_0(\resto{{\cal E}}{X^\st\times Y})$ has a distinguished extension.
\item  Consider a global solution $(H_t)_{t\in[0,\infty)}$ of the fiberwise Donaldson   evolution equation for the pair $({\cal E},K)$  
\begin{equation}\label{ev}
\frac{dH_t}{dt}=-2H_t (i\Lambda_Y [F_{H_t}^{Y,Y}]^0)\ ,\  H_t\in{\cal M}et^l({\cal E})\ ,\ H_0=K\ .
\end{equation}
The plurisubharmonic potential $\tilde\psi$ obtained in the proof of Theorem \ref{extension} satisfies the inequality
\begin{equation} \label{ineq}
\tilde \psi(x)\leq \varphi_K(x)+\lim_{t\to\infty} m(K,H_t)(x)=\varphi_K(x)+ \lim_{t\to\infty} M(K_x,(H_x)_t) \ \forall x\in X  .
\end{equation}
\item The completely pluripolar set $PP(pw_0({\cal E}))$ associated with the closed positive current $pw({\cal E})$ (i.e. the set where its local potentials take the value $-\infty$) satisfies
$$X\setminus X^\sst\subset PP(pw_0({\cal E}))\ .
$$
\end{enumerate}
\end{re}

There always exists a global solution $(H_t)_{t\in[0,\infty)}$ of the evolution equation (\ref{ev})   (see  (\cite{Do1}, \cite{DK} for $n=2$ and \cite{Si} Proposition 6.6 for the general case).  Moreover, if ${\cal E}_x$ is  stable, then $H_{t,x}$  converges weakly in $L^p_{2}$ to a projectively Hermitian-Einstein metric $H_\infty(x)\in {\cal M}et^{l_x}({\cal E}_x)$  as $t\to\infty$.  This second statement in the remark follows now by restricting $\tilde\psi$ to a small curve $C\subset X$ with $C\setminus\{x\}\subset X^\st$ and applying the mean value inequality for subharmonic functions. Indeed, note first that       
 $$\psi_t:=\varphi_{K}+m(K,H_t)\searrow \psi\hbox{ pointwise  on }  U\cap X^\st\hbox{ as }t\to\infty\ .$$
 
Identify  $U$  (which we suppose sufficiently small)  with an open set $V$ of $\C^{\dim(X)}$ such that $C\cap U$ is mapped onto an open set of a complex line. For sufficiently small $\xi\in\C$, the mean value inequality gives:
 $$\tilde\psi(x)\leq \frac{1}{2\pi} \int_0^{2\pi}  \psi (x+e^{i\theta}\xi)d\theta\leq \frac{1}{2\pi} \int_0^{2\pi}  \psi_t (x+e^{i\theta}\xi)d\theta\ .
 $$
Fix  $t>0$. Letting  $\xi$ tend to 0 and taking into account that $\psi_t$ is continuous, we get $\tilde\psi(x)\leq \psi_t(x)$. This holds for every $t$, so (\ref{ineq}) is proved.

 For the third statement, recall that, for a non-semistable bundle ${\cal F}$ on a compact K\"ahler manifold $(Y,g)$, the norm $\|\Lambda_Y(F_H^0)\|_{L^2}^2$ is bounded from below by a positive number when $H$ varies in the spaces of Hermitian metrics on ${\cal F}$ (see Proposition 5 and section 4 in \cite{Do1}). Using the identity  

$$\frac{d}{dt}Êm(K,H_t)=-2\| \Lambda_Y [F_{H_t}^{Y,Y}]^0\|^2 
$$
(see  \cite{Do1}), this implies  $\lim_{t\to\infty} M(K_x,H_{t,x})=-\infty$ when ${\cal E}_x$ is not semistable. It suffices now to use (\ref{ineq}).
\begin{re} \label{unitary} A similar extensibility result holds for the Petersson-Weil form $ pw(\resto{{\cal E}}{X^\st\times Y})$ associated with the holomorphic map $X^{\st}\to {\cal M}^\st(E)$ to the moduli space of all stable holomorphic structures on $E$. This result can be easily reduced to Theorem \ref{extension}, because $pw(\resto{{\cal E}}{X^\st\times Y})=pw_0(\resto{{\cal E}}{X^\st\times Y})+D^*(\omega_\Pic)$, where $D:X\to \Pic(Y)$ is given by $x\mapsto \det({\cal E}_x)$ and $\omega_\Pic$ is the (suitably normed) Petersson-Weil K\"ahler metric on $\Pic(Y)$.
\end{re}

\subsection{Remarks and open problems}
\label{remarks}

  Our extension theorem might look surprising: the Petersson-Weil metric extends (as a closed positive current) even over the locus where the bundles are not semistable. The result provides   strong evidence for the following conjectures:\\ \\
{\bf Conjecture 1:} A moduli space ${\cal M}^\st$ of stable bundles over a compact K\"ahler manifold admits  a   compactification $\bar {\cal M}^\st$ {\it which is a complex space} and is  obtained by adding to ${\cal M}^\st$ equivalence classes of sheaves satisfying a suitable semistability condition. This  compactification has the following natural property: For every generically stable family ${\cal E}$ parameterized by $X$, the induced holomorphic map $X^\st\to {\cal M}^\st$ extends to a meromorphic map $X\dasharrow  \bar {\cal M}^\st$.
\\ \\
This conjecture is not known even for K\"ahler surfaces. Jun Li's remarkable results \cite{Li} imply that the conjecture is true for algebraic surfaces. In this case one has two solutions of the problem: the Gieseker compactification and the Uhlenbeck compactification (which, surprisingly, also has the structure of a projective scheme). For K\"ahler surfaces the most natural approach is to use a Gieseker semistability condition obtained using formally the K\"ahler class instead of the polarization in the definition of the Hilbert polynomial. This approach has been suggested by Matei Toma. Another approach would be to try directly to put a complex space structure on the Uhlenbeck compactification of the moduli space. Results in this direction have been obtained by Matei Toma \cite{To} and Nicholas Buchdahl \cite{Bu}.  Note also that  Uhlenbeck-type compactness theorems have been obtained by Nakajima \cite{Na} and Tao-Tian \cite{Ti}, \cite{TT} for higher dimensional K\"ahler manifolds. This gives further evidence for the conjecture.  However  it is not known whether these compactifications have complex space structures.    
\\ \\
{\bf Conjecture 2:} The Petersson-Weil metric on the moduli space ${\cal M}^\st$ extends on a suitable complex geometric compactification  $\bar {\cal M}^\st$ endowing this compactification with the structure of a (in general singular) K\"ahler space.
\\ \\
This  only results we know  concerning this problem is due to Tyurin \cite{Ty}, who proved that, on K3 surfaces, the Petersson-Weil  metric extends smoothly on the Uhlenbeck compactification of certain moduli spaces of stable bundles.   Tyurin's result concerns only moduli spaces which contain no properly semistable bundles. This condition guarantees that the Uhlenbeck compactification contains no reduction.  On the other hand, even for algebraic surfaces, the conjecture seems to be very difficult. Indeed, in this case Jun Li's results  do give K\"ahler structures on both the Gieseker and the Uhlenbeck compactification, but these metric structures {\it cannot be easily compared to the Petersson-Weil metric}, and cannot be generalized to the non-algebraic K\"ahler framework. Note that, if the conjecture was true, the volume of the Petersson-Weil metric should be locally finite near a virtual instanton.  This problem can be studied even before proving Conjecture 1, but it  seems to be very difficult.
Our extension result (Theorem \ref{extension}) would follow easily from the two conjectures above.   Conversely, we believe that this result is an important first step towards solving these problems.
\begin{re} Conjecture 1 is definitely false for non-K\"ahlerian surfaces: Indeed, it is known that certain moduli spaces on class VII surfaces  with $b_2=1$ can be identified with open disks, and the only way to compactify these moduli spaces is to add the obvious boundaries of the disks consisting of split polystable bundles (see \cite{Sch}). Moreover, the author  showed that  the moduli space ${\cal M}^\pst(0,{\cal K})$ of polystable bundles ${\cal E}$ with $c_2({\cal E})=0$ and $\det({\cal E})={\cal K}$ on certain minimal class VII surfaces with $b_2=2$ is homeomorphic to $S^4$ \cite{Te3}.  
\end{re}

 Another class of interesting  problems suggested by our extension theorem is to relate the analytic  invariants of the closed positive current $pw({\cal E})$ to the complex geometry of our family ${\cal E}$. For instance:
 \begin{enumerate}
 \item Is the inclusion $X\setminus X^\sst\subset PP(pw_0({\cal E}))$ proved in Remark \ref{interesting} an equality?
 \item \label{computelelong} Compute the Lelong numbers of $\nu(pw_0({\cal E}),x)$, $x\in X\setminus X^\st$  in terms of complex  geometric properties of ${\cal E}$.   
 \item \label{ideal} Describe the Nadel multiplier ideal sheaf ${\cal I}(pw_0({\cal E}))$ of the current in terms of complex  geometric properties of ${\cal E}$.  
  The complex subspace defined by this ideal sheaf might be called the {\it non-stability complex subspace}  of the family.
 \end{enumerate}
 We believe that for solving problems (\ref{computelelong}) and (\ref{ideal}) an important role is played by the Harder-Narasimhan filtrations of the restrictions  $\resto{{\cal E}}{Y^{(m)}_x}$    to the higher order  infinitesimal neighborhoods $[Y_x]^{(m)}$ of the fibers $Y_x$ over points $x\in X\setminus X^\sst$.

 \section{Families of bundles over a K\"ahler manifold parameterized by a non-K\"ahlerian surface}
 \label{incompatsection}
 
The purpose of this section is to prove Theorem \ref{incompatibility} stated in the introduction. This result  and Corollary \ref{goal} play  an important role in our work  about existence of curves on class  VII surfaces with $b_2=2$.
\\
\\
\pf (of Theorem \ref{incompatibility}). Using Remark \ref{unitary} consider the Petersson-Weil current $pw({\cal E})$ of the family, and let $R$ be its residual part with respect to Siu's decomposition formula (see  for instance  2.18 \cite{De2}). $R$ and  $pw({\cal E})$ coincide on  $X^\st$, where both are smooth. By Corollary \ref{dan} in the Appendix $R$ cannot be strictly positive and smooth on a non-empty open set. Therefore the smooth form $pw(\resto{\cal E}{X^\st\times Y})$ is degenerate, which implies that the derivative $f_*$ of the classifying map $f:X^\st\to{\cal M}^\st(E)$ has rank at most 1 at every point. This proves (1).

Suppose now that   $f$ has generically rank  1, and let $X^1\subset X^\st$ the non-empty Zariski open subset of points where the rank  is exactly 1. The image  $f(X^{1})\subset {\cal M}^\st(E)$ contains infinitely many isomorphism classes. For every $[{\cal F}]\in f(X^{1})$, consider the Brill-Noether locus
$$X_{\cal F}:=\{x\in X|\ H^0(p_Y^*({\cal F})^\vee\otimes{\cal E})\ne 0\}\subset X\ . 
$$
For any $[{\cal F}] \in  f(X^{1})$ one has: 
\begin{enumerate}
\item $X_{\cal F}$ is an analytic subset of $X$,
\item \label{curve} $X_{\cal F}\cap X^1$ is a non-empty curve.
\item For $[{\cal F}]\ne [{\cal F}']\in  f(X^{1})$ one has $X_{\cal F}\cap X_{{\cal F}'}\cap X^\st=\emptyset$.
\end{enumerate}

Let $C_{\cal F}$ be a 1-dimensional irreducible component of $X_{\cal F}$ which intersects $X^1$. Such a component exists by (\ref{curve}). For   $[{\cal F}]\ne [{\cal F}']$ one must have $C_{\cal F}\ne C_{{\cal F}'}$ (because they are disjoint on $X^\st$).
\qed

Using Theorem \ref{incompatibility} we can finally prove Corollary \ref{goal} stated in the introduction:\\

\pf (of Corollary \ref{goal}) If the family was generically stable, the corresponding classifying map $f:X^\st\to {\cal M}^\st$ must be constant. Let $[{\cal F}_0]\in {\cal M}^\st$ be this constant. Then  the sheaf ${\cal U}:=[p_X]_*(p_Y^*({\cal F}_0^\vee)\otimes {\cal E})$ has rank 1, because it is a line bundle on $X^\st$ by Grauert's local-triviality Theorem. One obtains a tautological morphism
$$[p_X]^*({\cal U})\otimes p_Y^*({\cal F}_0) \map {\cal E}\ ,
$$
which is a bundle isomorphism on $X^\st\times Y$. Restricting this tautological morphism to fibers $X\times\{y\}$, we get morphisms 
$${\cal U}\otimes {\cal O}^{\oplus 2}_X\simeq {\cal U}\otimes {\cal F}_0(y)\to {\cal E}_y$$
which are bundle isomorphisms on $X^\st$. We define ${\cal T_0}$ to be the reflexivization of the left hand sheaf. \\
 
Suppose now that $X^\st=\emptyset$.  In this case we have a family of filtrable bundles, so by Remark \ref{oldremark} the set of jumps $J:=\{\degmax({\cal E}_x)|\ x\in X\}$ is finite. Let $d_0$ be the minimal element. One has $X_{\geq d_0}=X$, and the  set $V:=X\setminus X_{>d_0}$  is Zariski open and non-empty. Note that
$$V=\{x\in X|\ \degmax({\cal E}_x)=d_0\}\ .
$$
Consider the following diagram
$$\begin{array}{ccc}
BN_X({\cal E})_{\geq d_0}&\map& \Pic(Y)\\
\\
q_X\downarrow\phantom{p_X}&& \\
\\
X&&
\end{array}
$$
with compact left hand terms. By Lemma \ref{nonstable} below, the vertical map is one to one or two to one   above $V$, so the diagram defines a meromorphic map $\varphi:X\dasharrow \Pic(Y)$ (or $\varphi:X\dasharrow S^2(\Pic(Y))$) which is regular on $V$ (or an smaller dense Zariski open set). But $\Pic(Y)$ is K\"ahler, so this map induces a smooth closed positive (1,1)-form on $V$ which extends as closed positive current on $X$. The same arguments as in Theorem \ref{incompatibility} show that $\varphi$ must be constant on $V$. Let ${\cal L}_0$ (or $\{{\cal L}_0,{\cal L}_0'\}$) be this constant. By Lemma \ref{nonstable} and Grauert's local triviality theorem, we see that the sheaf
$${\cal T}:=[p_X]_*(p_Y^*({\cal L}_0^\vee)\otimes {\cal E})$$
has rank 1 or 2 and is locally free on a dense, Zariski open subset  $U\subset V$ (the subset of points $v\in V$ where $h^0({\cal L}_0^\vee\otimes {\cal E}_v)$ is minimal). We get again a tautological map
$$p_X^*({\cal T})\otimes p_Y^*({\cal L}_0)\map {\cal E}
$$
which is a bundle embedding on $U\times Y$. Restricting this morphism to fibers $X\times\{y\}$, and taking reflexivizations, we get the result.
\qed
\begin{lm} \label{nonstable} Let ${\cal E}$ be a non-stable rank 2 bundle on a compact Gauduchon manifold $(Y,g)$.  Put 
$$BN^{\max}({\cal E}):=\{{\cal L}\in\Pic(Y)|\ \deg({\cal L})=\degmax({\cal E}),\ H^0({\cal L}^\vee\otimes{\cal E})\ne 0\}\ .$$
1. One of the following holds
\begin{enumerate}
\item ${\cal E}$ is non-semistable. In this case $BN^{\max}({\cal E})=\{{\cal L}_{\max}\}$, where ${\cal L}_{\max}$ is the unique maximal destabilizing line bundle of ${\cal E}$.
\item ${\cal E}$ is semistable and can be written as a nontrivial extension 
\begin{equation}\label{exactseq}0\map {\cal L}\map {\cal E}\map {\cal M}\otimes {\cal I}_Z\map 0\ ,
\end{equation}
where $Z\subset Y$ is a (possibly empty) 2-codimensional locally complete  intersection  and $\deg({\cal L})=\deg({\cal M})=\frac{1}{2}\deg({\cal E})$. In this case 
$BN^{\max}({\cal E})=\{{\cal L}\}$.
\item ${\cal E}$ is semistable and is isomorphic to the direct sum of two non-isomorphic line bundles ${\cal L}_1$, ${\cal L}_2$. In this case $BN^{\max}({\cal E})=\{{\cal L}_1,{\cal L}_2\}$.
\item ${\cal E}$ is semistable and is isomorphic to ${\cal L}\oplus  {\cal L}$ for a line bundle ${\cal L}$. In this case $BN^{\max}({\cal E})=\{{\cal L}\}$.
\end{enumerate}
\vspace{1mm} 
2.  Let   ${\cal L}_0\in\Pic(Y)$ with $\deg({\cal L}_0)=\degmax({\cal E})$. Then $h^0({\cal L}_0^\vee\otimes {\cal E})\leq 2$. Equality holds  if only if ${\cal E}$ is the split polystable bundle ${\cal L}_0\oplus{\cal L}_0$. 

\end{lm}
\pf Recall  that reflexive rank 1 sheaves are invertible. Therefore any rank 1 subsheaf of ${\cal E}$ is contained in the kernel of an exact sequence of the form (\ref{exactseq}).  Since ${\cal E}$ was assumed to be non-stable, there exists a destabilizing short exact sequence of this form   with $\deg({\cal L})=\degmax({\cal E})\geq\deg({\cal M})$. For any line bundle ${\cal L}_0\in\Pic(Y)$ with $\deg({\cal L}_0)=\degmax({\cal E})$ one has $\deg({\cal L}_0)\geq \deg({\cal M})$, so $\Hom({\cal L}_0,{\cal M}\otimes {\cal I}_Z)=0$ except when $Z=\emptyset$, ${\cal L}_0\simeq {\cal M}$ and the extension (\ref{exactseq}) splits (i.e. when we are in   case (3) or (4)). In   case (1) or (2), we will have $\Hom({\cal L}_0,{\cal M}\otimes {\cal I}_Z)=0$, so $\Hom({\cal L}_0,{\cal E})=\Hom({\cal L}_0,{\cal L})$, which is either 1-dimensional   (when ${\cal L}_0\simeq  {\cal L}$), or vanishes.
\qed

 \section{Appendix. A self-intersection inequality} \label{appendix}
 
 In this section we prove a technical result concerning the self-intersection of a closed positive  current on arbitrary (possibly non-K\"ahlerian) surfaces. The author is indebted to Dan Popovici, S\'ebastien Boucksom and Tien-Cuong Dinh for interesting and useful discussion on the subject.  The idea to use Boucksom's recent results on approximation of closed positive currents by currents with analytic singularities is due to Dan Popovici. Our goal  is the following:
   \begin{thry} \label{self} Let $T$ be a closed positive (1,1)-current  on a compact complex surface $X$ with the property that the associated analytic sets $E_c(T)$ are all 0-dimensional (or empty) for every $c>0$. Let $M\subset X$ be the countable set 
 $$M:=\{x\in X|\ \nu(T,x)>0\}\ .
 $$
 Then the square $[T]^2$ of the cohomology class $[T]$ of $T$ is represented by a positive current $\Theta$ satisfying the inequality
 $$\Theta\geq \sum_{m\in M} \nu(T,m)^2[m]+T_{\rm ac}\wedge T_{\rm ac}\ .
 $$
 \end{thry}

 This result is is a special case of Demailly's self-intersection  theorem (see Theorem 9.5 in \cite{De2} or Theorem 1.7, Corollary 7.6  in \cite{De1}). The proof given in these articles requires K\"ahlerianity, but the special case stated here does not. For completeness we insert a short proof of Theorem  \ref{self}.
  
 Let $U$ be an arbitrary complex manifold (not necessary compact).  We recall that for an (almost)plurisubharmonic function $u$ which is locally bounded on the complement of  an analytic 0-dimensional set  $A\subset U$, and a closed (almost) positive current (1,1) $T$ on $U$ one can define in a coherent way the wedge product $dd^c u\wedge T$ (see \cite{De2} Proposition 2.3). We will call this operation {\it generalized wedge product}. On the other hand, for an arbitrary {\it continuous} form $\eta$, the wedge product $\eta\wedge T$ is also defined, because closed almost positive currents have measure coefficients. The two wedge operations are defined  in different ways. However  one has

 \begin{lm}\label{obvious}  Let $\eta$, $\mu$ be continuous $(1,1)$-forms, $T$ a closed almost positive current, and $u$ an almost  plurisubharmonic function which is locally bounded on the complement of  an analytic 0-dimensional set. The current 
 $$dd^c u\wedge T+\eta\wedge T+ dd^c u\wedge \mu +\eta\wedge\mu$$
\begin{enumerate}
 \item depends only on $dd^c u+\eta$ and $T+\mu$. 
 \item Is positive when $dd^c u+\eta$ and  $T+\mu$ are positive.
\end{enumerate}
  \end{lm}
  \pf  Since the problem is local, we may suppose that $U\subset\C^n$. Let $\rho_\varepsilon$, $\varepsilon\in]0,1[$ be a system of smoothing kernels. Suppose $dd^c u+\eta=dd^c u'+\eta'$, $T+\mu=T'+\mu'$ with $\eta$, $\eta'$, $\mu$, $\mu'$ continuous.  Using the weak continuity property of the generalized product, it follows that
on every relatively compact open set $V\subset U$ it holds
  $$dd^c u\wedge T+\eta\wedge T+ dd^c u\wedge \mu +\eta\wedge\mu=\lim_{\varepsilon\searrow 0}(dd^c u_\varepsilon\wedge T+\eta_\varepsilon\wedge T+dd^c u_\varepsilon\wedge\mu+\eta_\varepsilon\wedge\mu)=$$
  $$\lim_{\varepsilon\searrow 0}(dd^c u_\varepsilon+\eta_\varepsilon)\wedge(T+\mu)=\lim_{\varepsilon\searrow 0}(dd^c u'_\varepsilon+\eta'_\varepsilon)\wedge (T'+\mu')=$$
  $$=dd^c u'\wedge T'+\eta'\wedge T'+ dd^c u'\wedge \mu' +\eta'\wedge\mu'\ .
  $$
This proves the first statement. Suppose now that $dd^c u+\eta\geq 0$, $T+\mu\geq 0$. Then (since positivity is preserved by convolution) one obtains $dd^c u_\varepsilon+\eta_\varepsilon\geq 0$ on an open subset  $U_\varepsilon$, where the system $(U_\varepsilon)_{\varepsilon>0}$ is an exhaustion of $U$.  On $U_\varepsilon$ we can write
$$dd^c u_{\varepsilon}\wedge T+\eta_\varepsilon\wedge T+dd^c u_\varepsilon\wedge \mu+\eta_\varepsilon\wedge\mu=(dd^c u_{\varepsilon}+\eta_{\varepsilon})\wedge (T+\mu)\geq 0\ .$$
One a fixed relatively compact open subset $V\subset U$ one has $\eta_\varepsilon\to \eta$ uniformly, (because $\eta$ is continuous) and $dd^c u_\varepsilon\to dd^c u$ in the space of  measures.  Therefore  the left hand side converges weakly to  $dd^c u\wedge T+\eta\wedge T+ dd^c u\wedge \mu +\eta\wedge\mu$.
\qed

 Lemma \ref{obvious} has a global version.
 \begin{lm}\label{easy} Let $X$ be a complex manifold and  $T$ a closed almost positive $(1,1)$-current    and  $\eta$, $\mu$  continuous (1,1)-forms on $X$. For any closed almost positive $(1,1)$-current  $S$  admitting    locally almost plurisubharmonic   potentials which are locally bounded on the complement of a 0-dimensional analytic set, the current $S\wedge T+\eta\wedge T+S\wedge\mu+\eta\wedge\mu$ is well defined, depends only on   $S+\eta$ and $T+\mu$, and is positive when  $S+\eta$ and $T+\mu$ are positive. Moreover, when $X$ is a compact surface, then
 $$\langle S\wedge T+\eta\wedge T+S\wedge\mu+\eta\wedge\mu, 1\rangle=\langle [S]\cup[T], [X]\rangle +\langle T,\eta\rangle+\langle S,\mu\rangle+\int_X \eta\wedge \mu  \ ,
 $$
 where $[S]$, $[T]$ denote the de Rham cohomology classes of $S$ and $T$.
 \end{lm}
The integral formula in the lemma follows from Corollary 9.2 in \cite{De2}.
\\ \\
 \pf   (of Theorem \ref{self})     If  $T$  is a closed positive current with   analytic singularities, it is well-known that its Siu decomposition $T=\sum_n\nu_n [C_n]+R$ coincides with its Lebesgue decomposition $T=T_s+T_{ac}$ \cite{B1}. In our case, since all the  analytic sets $E_c$ are 0-dimensional, we have  $T = R = T_{ac}$, and therefore $T$ is absolutely continuous.  By Proposition 2.3 in \cite{De2}, the square $\Theta=T^2$ is a well-defined  closed positive (2,2)-current representing the cohomology class $[T]^2$. Although $T$ is absolutely continuous, its generalized wedge square cannot be computed pointwise.
 
 Let  $M\subset X$ the finite set of singularities of $T$. On $X\setminus M$ the current $T$ (hence also $\Theta=T^2$) is smooth.  This shows that $\Theta_{ac}=T_{ac}\wedge T_{ac}$. Moreover, the multiplicative property of the Lelong numbers (Proposition 2.16 \cite{De2}) shows that $\nu(\Theta,m)\geq \nu(T,m)^2$ for every $m\in M$. Therefore $\Theta\geq \sum_{m\in M} \nu(T,m)^2 [m]+ T_{ac}\wedge T_{ac}$.

Now suppose that $T$ has arbitrary singularities. Choose a Gauduchon metric $\omega$ on $X$. Such a metric always exists. Demailly's regularization theorem (\cite{De1}, Prop. 3.7., p. 380) yields a sequence of almost  positive  closed (1,1)-currents $T_k$ with analytic singularities in the $dd^c$-cohomology class of $T$ satisfying: 
\begin{enumerate}
\item  $T_k \to T$ in the weak topology of currents as $m\to\infty$; 
\item $T_k \geq -\varepsilon_k\omega$, for a sequence $\varepsilon_k$ decreasing to 0; 
\item  \label{lelongsmaller} $\nu (T_k, x) \leq \nu (T,x)$ for every $x \in X$ and the function $\nu(T_k,\cdot)$ converges uniformly to $\nu(T,\cdot)$; 
\item  \label{bouck} $(T_k)_{ac}(x)$ converges to $T_{ac}(x)$ for almost every $x\in X$. 
\end{enumerate}
Property (\ref{bouck}) is not part of Demailly's original theorem. It is a later addition 
of Boucksom (Theorem 2.6, Corollary 2.7 \cite{B2}).   Using (\ref{lelongsmaller}) we see that $E_c(T_k)\subset E_c(T)$ (which is 0-dimensional) for every $c>0$,  so the singular set of $T_k$ is a finite set $M_k$, and $T_k$ is smooth on $X\setminus M_k$. By  Lemma \ref{easy} the square  $\Theta_k:=(T_k+\varepsilon_k\omega)^2$ is well defined and positive, and its total mass is
$$\langle \Theta_k, 1\rangle =[T_k]^2+2 \varepsilon_k \langle T_k,\omega \rangle+\varepsilon_k^2 \int_X \omega ^2
$$

We have obviously $[T_k]^2=[T]^2\in H^4(X,\R)$ and (since $\omega$ is Gauduchon) $ \langle T_k,\omega \rangle=\langle T,\omega\rangle$. The total mass of  $\Theta_k$ is bounded so (passing to a subsequence if necessary), we get a weak limit $\Theta:=\lim_{k\to \infty}\Theta_k$.  Since $\Theta_k$ and $T_k$ are both smooth on $X\setminus M_k$ we get $(\Theta_k)_{ac}=(T_k)_{ac}\wedge (T_k)_{ac}+ 2\varepsilon_k  (T_k)_{ac}\wedge \omega+ \varepsilon_k^2\omega^2$, which converges to $T_{ac}\wedge T_{ac}$ almost everywhere by (\ref{bouck}) in Demailly-Boucksom's regularization theorem. Let $f$ be a non-negative continuous function on $X$. By Fatou's lemma we get
$$\int\limits_X f T_{ac}\wedge T_{ac}=\int\limits_X  \liminf_{k\to\infty} f (\Theta_k)_{ac}\leq \liminf_{k\to\infty}  \int\limits_X f (\Theta_k)_{ac}\leq  \liminf_{k\to\infty}  \int\limits_X f (\Theta_k) =\langle \Theta, f\rangle\ .
$$
Taking $f=1$, we see that the non-negative measurable 4-form $T_{ac}\wedge T_{ac}$ belongs to $L^1$. Letting $f$ vary in the space of non-negative continuous functions, we get   $\Theta\geq T_{ac}\wedge T_{ac}\in L^1$, which shows that
\begin{equation}\label{acpart}
\Theta_{ac}\geq T_{ac}\wedge T_{ac}\ .
\end{equation}
The Lelong numbers of $\Theta$ can be estimated  as follows. First note that
$$\Theta=\lim_{k\to\infty}T_k \wedge  T_k + 2\varepsilon_k   T_k \wedge \omega+ \varepsilon_k^2\omega^2= \lim_{k\to\infty}T_k \wedge  T_k\ ,
$$
because the term $\varepsilon_k   T_k \wedge \omega$ converges weakly to 0. Let $x\in X$ and $\omega_0$ be a closed K\"ahler form defined on an open  neighborhood $U$ of $X$ such that $\omega_0\geq \resto{\omega}{U}$. One can write $\resto{\Theta}{U}=\lim_{k\to\infty}\resto{T_k}{U}^2= \lim_{k\to\infty} (\resto{T_k}{U}+\varepsilon_k\omega_0)^2$. Therefore, using the semicontinuity and multiplicative properties of  the Lelong numbers associated to closed   positive  currents \cite{De2}, we get
\begin{equation}\label{spart}
\nu(\Theta,x)\geq \limsup_{k\to \infty} \nu(T_k,x)^2=\nu(T,x)^2\ .
\end{equation}

The statement follows now from (\ref{acpart}) and (\ref{spart}).
\qed

\begin{co}\label{dan}
In the conditions of Theorem \ref{self} one has
$[T]^2\geq \int_X T_{ac}\wedge T_{ac}\ .$
 If, moreover,    $T$ is continuous and strictly positive on non-empty open set, then $[T]^2>0$. In particular the surface must be K\"ahlerian.

\end{co}

\vspace{1mm} 
{\small
Author's address: \vspace{2mm}\\
Andrei Teleman, LATP, CMI,   Universit\'e de Provence,  39  Rue F.
Joliot-Curie, 13453 Marseille Cedex 13, France,  e-mail:
teleman@cmi.univ-mrs.fr. }

\end{document}